\magnification=\magstep1
\input amstex
\documentstyle{amsppt}
\topmatter \title Toda and KdV \endtitle
\author D. Gieseker \endauthor
\affil Department of Mathematics, U.C.L.A. \endaffil
\address Los Angeles, CA 90095\endaddress
\email dag@math.ucla.edu \endemail
\thanks I wish to thank Pantelis Damianou, Ed Frenkel, and Ali Kisisel to helpful conversations. This work was partially supported by the N.S.F.\endthanks
\abstract The main object of this paper is to produce a deformation of the KdV hierarchy of partial differential equations. We construct this deformation by taking a certain limit of the Toda hierarchy. This construction also provides a deformation of the Virasoro algebra.\endabstract
\endtopmatter
\document
\head 1. Introduction   \endhead
Our aim in this paper is to produce a deformation of the KdV hierarchy whose existence was conjectured in \cite {G}. To describe KdV algebraically following Gelfand and Dickey \cite {D}, let $$R_0 =\bold C[w^{(0)},w^{(1)}\ldots]$$ be a polynomial ring in infinitely many variables.  Introduce a $\bold C$ derivation $\partial$ on $R_0$ by $$\partial w^{(k)}= w^{(k+1)}.$$ An element of $R_0$ is intended to represent an abstract differential operator in one variable.  If $f $ is a $\Cal C^{\infty}$ function on $\bold R$, then  define $$P(f) = P(f,\frac{d f}{d x} ,\ldots),$$ i.e.  substitute $f$ for $w^{(0)}$, $\frac{d f}{d x}$ for $w^{(1)}$, etc.  

To describe translationally invariant PDE's algebraically,  consider $\bold C$ derivations $D$ of $R_0$ which commute with $\partial.$  The set of such $D$ is naturally just $R_0$ under the correspondence $D \to D(w^{(0)}).$ So $R_0$ inherits the structure of a Lie algebra, since the commutator of two derivations is a derivation. Let $$K_1=w^{(3)}+w^{(1)}w^{(0)},$$ called the KdV element of $R_0.$   One of the main results of KdV theory is that $K_1$ lies in a large abelian subalgebra of $R_0$. In fact, there is a sequence of elements $K_n$ of $R_0$ so that $K_{n+1}$ is not in the $\partial$ invariant subring of $R_0$ generated by the lower $K_n \ldots K_1$ and their derivatives and all the $K_k$ commute.  These $K_k$ are called the KdV hierarchy.

The main object of this paper will be to produce interesting deformations of the KdV hierarchy. That is, we seek to produce mutually commuting $L_k\in R_0[[\epsilon]]$ which become the $K_k$ when we set $\epsilon =0,$ where   $R_0[[\epsilon]]$ is the formal power series ring in $\epsilon.$ 

While the ring $R_0$ captures much of the algebraic structure of KdV, sometimes one considers solutions  or representations of KdV.  Suppose that $\Cal M$ is an analytic manifold and let $P_0 =\partial , P_1, \ldots $ be derivations of $R_0$ commuting with $\partial.$  Let $\chi_0,\chi_1,\ldots, $ be vector fields on 
$\Cal M$ and let $f$ be a function on $\Cal M.$  We say that $f,\chi_0,\chi_1\ldots$ form a representation of $P_0,P_1,\ldots$ if
$$\chi_k(f)=P_k(f,\chi(f),\chi^2(f) \ldots),$$ where we regard the vector field $\chi_k$ as a derivation on functions on $\Cal M$ and $\chi =\chi_0.$  

The main example of representations of the $K_k$ is the following: 
Let $X$ be a hyperelliptic curve and let $Q$ be a Weirstrass point. Let $$\vartheta:H^1(\Cal O_X) \to \bold C$$ be the theta function.  Then we can find translationally invariant vector fields, $\chi_0,\chi_1 \ldots$ so that if 
$$f = \chi_0^2(\log \vartheta),$$ then $f,\chi_0\ldots$ form a representation of the KdV hierarchy.

Our aim is to develop a difference version of KdV hierarchy to obtain a deformation of the KdV hierarchy. Here the idea is basically to discretize a differential equation.  The heuristic motivation for these discretizations in given in \cite{G}. To be rigorous  we need a method to describe  difference equations.  We consider the ring $$S_1=\bold C[\ldots X_{-1},X_0,X_1\ldots;\ldots Y_{-1},Y_0,Y_1\ldots]$$ and let  $$S_2=\bold C[\ldots a_{-1},a_0,a_1\ldots;\ldots b_{-1},b_0,b_1\ldots].$$ Let $T:S_2 \to S_2$ be the $\bold C$ algebra homomorphism defined by $T(a_n)=a_{n+1}$ and $T(b_n)=b_{n+1}.$ Now given $P,Q \in S_1$, we can define a derivation $D_{P,Q}:S_2\to S_2$ by 
$$D_{P,Q}(a_n)=P(\ldots,a_{n-1},\hat a_n , a_{n+1}\ldots ;\ldots,b_{n-1},\hat b_n , b_{n+1}\ldots)$$ and 
$$D_{P,Q}(b_n)=Q(\ldots,a_{n-1},\hat a_n , a_{n+1}\ldots ;\ldots,b_{n-1},\hat b_n , b_{n+1}\ldots),$$
where the  $\hat{\phantom{x}}$ indicates that $a_n$ should be substituted for $X_0$  and $b_n$ should be substituted  for $Y_0.$ This construction gives all the derivations of  $S_2$ commuting with $T$ and so introduces a Lie algebra structure on $S_1 \oplus S_1.$ The interesting example is the Toda equations:
$$T_1=(P_1,Q_1)=(Y_{-1}-Y_0,Y_0(X_0-X_1)).$$  The main theorem here due to Toda, Flaschka and many others is that $T_1$ lies in an unexpectedly large Abelian sub-algebra of the Lie algebra $S_1 \oplus S_1.$ In fact, there is a whole sequence of mutually commuting $T_k \in S_1 \oplus S_1.$ 

We can also describe solutions of the Toda hierarchy using algebraic geometry following van Moerbeke. Let $N$ be a positive integer. Let $\Cal C$ be the space
of all complex valued functions on $\bold Z$.  Let $T:\Cal C \to \Cal C$ be
translation by $N$, $T(f)(n)=f(n+N)$.  Let $\Cal C_N$ be the set of
translation invariant functions: $T(f)=f.$ Given $A$ and $B$ in
$\Cal C_N$, we define $$L_{(A,B)}:\Cal C \to \Cal C$$ by the formula 
$$L_{(A,B)}(\psi)(n)=\psi(n+1)+A(n)\psi(n)+B(n)\psi(n-1).$$ Thus
$L_{(A,B)}(\psi)$ is a second order linear difference operator.    By
definition of $\Cal C_N$, the operators $L_{(A,B)}$ and $T$
commute, so we can reasonably look for common eigenfunctions of these
two operators.  If you think of $L_{(A,B)}$ as a discrete analogue
of a Schr\"odinger operator, this amounts to finding the energy levels
with a given quasi-momentum of a particle traveling through a
periodic potential, a problem frequently encountered in solid state
physics  \cite{AM}.  We then define the Bloch spectrum $\Cal B_{(A,B)}$
of $L_{(A,B)}$ to be the set of $(\lambda,\alpha)\in \bold C \times
\bold C^*$ so that there is a non-zero function $\psi$ with
$L_{(A,B)}(\psi)=\lambda \psi$ and $T(\psi)=\alpha \psi.$ By
projecting to the $\lambda$ axis, it is easy to see that
$\Cal B_{(A,B)}$ is a hyperelliptic curve, possibly singular.  Indeed,
there are in general two
values of $\alpha$ associated to any fixed  $\lambda.$ $\Cal B_{(A,B)}$ can
be compactified to a curve $\bar
\Cal B_{(A,B)}$ by adding two points $P$ and $Q$ over $\lambda=\infty.$ It turns out that the divisor $N(P-Q)$ is linearly equivalent to zero. It also turns
out that $\Cal B_{(A,B)}$ does not determine $(A,B)$.  There are
interesting ways of moving $(A,B)$ so that $\Cal B_{(A,B)}$ remains
fixed.  Such a deformation of $(A,B)$ keeping the Bloch spectrum fixed
is called an isospectral deformation.  The set of all $(A',B')$
isospectral to $(A,B)$ turns out to be isomorphic to the Jacobian of  $\bar
\Cal B_{(A,B)}$ in a birational sense for generic $A$ and $B$. In particular, any linear flow on the Jacobian becomes a non-linear flow on $\Cal C_N \times \Cal C_N$, which in turn is a linear combination of Toda flows.  For example, there is a linear flow on the Jacobian so that if $(A_t,B_t)$ indicates the flow of $(A_0,B_0)$ after time $t$, then 
$$\frac{dA_t(k)}{dt} = B_t(k-1)-B_t(k)$$ and 

$$\frac{dB_t(k)}{dt} = B_t(k) (A_t(k)-B_t(k+1)).$$ Conversely, given a hyperelliptic curve $C$ of genus $g$ and two points $P$ and $Q$ on $C$ and a suitably generic line bundle $\Cal L$ of degree $g$ so that $N(P-Q)$ is linearly equivalent to zero, we can define $A_{\Cal L} $ and $B_{\Cal L}$ in $\Cal C_N$ so that linear flows on the Jacobian become Toda flows on $\Cal C_N \times \Cal C_N.$ (See \cite {MM} for instance.)

\cite{G} developed a method for taking the limit of Toda equations.  Our purpose in this paper is to  develop an algebraic framework for these limits and to prove the main conjecture of \cite{G}. Here is a formalism to allow us to make sense of taking a limit of Toda equations. Let 
$$R =\bold C[v^{(0)},v^{(1)}\ldots; w^{(0)},w^{(1)}\ldots].$$ Again we introduce a derivation $\partial$ by $\partial v^{(k)}=v^{(k+1)}$ and by considering $R_0$ as a subring of $R.$  We think of elements $P \in R$ as being differential expressions in two functions $f$ and $g$:
$$P(f,g) =P(f,\frac{df}{dx},\ldots;g,\frac{dg}{dx}\ldots).$$ For $k\in \bold Z$, let $E_k:R[[\epsilon]] \to R[[\epsilon]]$ be defined by 
$$E_k =\exp(k\epsilon \partial)$$ as a formal power series in $\epsilon$. Now suppose we have $(P_1,P_2)\in S_1 \oplus S_1.$ We can then define a derivation $\Cal D'_{P_1,P_2}:R[[\epsilon]] \to R[[\epsilon]]$ commuting with $\partial$ and continuous in $\epsilon$ topology 
by $$\split \Cal D'_{P_1,P_2}(v^{(0)})= \\& P_1(\ldots
E_{-1}(v^{(0)}),\hat E_{0}(v^{(0)}),E_{1}(v^{(0)}),\ldots;\\& \ldots,
E_{-1}(w^{(0)}),\hat E_{0}(w^{(0)}),E_{1}(w^{(0)}),\ldots)\endsplit$$
and $$\split \Cal D'_{P_1,P_2}(w^{(0)})= \\& P_2(\ldots
E_{-1}(v^{(0)}),\hat E_{0}(v^{(0)}),E_{1}(v^{(0)}),\ldots;\\& \ldots,
E_{-1}(w^{(0)}),\hat E_{0}(w^{(0)}),E_{1}(w^{(0)}),\ldots)\endsplit.$$
Define a $\bold C$ algebra endomorphism of $R[[\epsilon]]$ commuting
with $\partial$ by $$\Phi(v^{(0)})=-2+\epsilon^2 v^{(0)}$$ and
$$\Phi(w^{(0)})=1+\epsilon^2 w^{(0)}.$$ $\Phi$ does not have an
inverse on $R[[\epsilon]]$, but does have one on $R((\epsilon))$, the ring of
Laurent series in $\epsilon$ which contain only finitely many negative
powers of $\epsilon$ and finally define 
$$\Cal D_{P_1,P_2} = \Phi \Cal D'_{P_1,P_2} \Phi^{-1}.$$ In particular, it turns out that we get a series of mutually commuting derivations $\Cal D_{T_k}$ of $R[[\epsilon]]$ coming from the Toda hierarchy. 

Suppose that $P \in R_0[[\epsilon]]$ and consider the element $W=v^{(0)}-P$ and let $\Cal I_P$ be the closure of the ideal of $R[[\epsilon]]$ generated by $W$, $\partial W$, $\partial^2 W$, etc. Notice that $R[[\epsilon]]/\Cal I_P$ is naturally isomorphic to $R_0[[\epsilon]].$ The following is one of our  main results:
\proclaim{Theorem 1.1} There is a $P$ so that $\Cal I_P$ is invariant under all the $\Cal D_{T_k}.$ Thus the $\Cal D_{T_k}$ induce derivations $\bold  D_{T_k}$ of $R_0[[\epsilon]].$ Suitable linear combinations of the $\bold  D_{T_k}$ over $\bold C((\epsilon))$ are a deformation of the KdV hierarchy. \endproclaim

We construct this $P$ recursively in powers of $\epsilon$. Suppose we have found a $P$ which works to order $\epsilon^n$ and we write
$P'=P+\epsilon^n P_1.$ For  $\Cal D_{T_1}$ to preserve $\Cal I_{P'}$ to order $\epsilon^{n+1}$ it turns out that there is an element $Q_n$ computed in terms of the original $P$ so that $\partial P_1 =Q_n.$  But it is not at all obvious why $Q_n$ should be the total derivative of 
anything. Thus at each stage of constructing $P$, we meet a highly non-trivial obstruction. We will show that there are lots of functions $g$ so that 
$$\int_z^{z+1} Q_n(g)=0\tag 1.1.1$$ for generic $z\in \bold C.$  If 1.1.1 is true for generic enough $g$, then $Q_n$ is a derivative.   

We define a generalized idea of representations. Suppose that we have an element  $P \in R$.   Let
$\Cal M$ be an analytic manifold and let $\chi$ be an analytic vector
field on $\Cal M$ and let $f$ and $g$ be two meromorphic functions on
$\Cal M$.  We define $P_\chi(f,g)$ to be $P(f,\chi f, \chi^2 f\ldots g, \chi g,\chi^2 g \ldots).$ Next, suppose  we have a function $h$ on $\Cal M.$ For any given $N\in \bold Z^{+}$, 
we can extend this definition to $P \in
R[[\epsilon]]$ by 
$$P_{\chi,N} (\sum_{n=0}^\infty P_n\epsilon^n)(f,g) =\sum_{n=0}^N
h^n P_{n,\chi}(f,g).$$ If $h_1$ and $h_2$ are two functions on $\Cal M$,
we say 
$$h_1 \equiv h_2 \mod h^N$$ if $(h_1 -h_2)/h^N$ is analytic at all the points $P$
where both $h$ is equal to zero at $P$ and both $h_1$ and $h_2$ are analytic at
$P$. We will use a similar terminology for vector fields on $\Cal M$.

Suppose we are given  derivations $D_1 \ldots D_n$ of
$R[[\epsilon]]$ and we are given vector fields $\chi,\chi_1 \ldots \chi_n$ on a manifold $\Cal M.$
Suppose we are also given a function $h$ on $\Cal M$ which is killed
by $\chi$ and all the $\chi_j.$ We further suppose that the $f$ and
$g$ do not have poles along the set $\{h=0\}$.

\proclaim{Definition 1.2} In the above situation, we say that 
$(f,g,h;\chi,\chi_1 \ldots \chi_n)$  form a 
representation 
of $D_1 \ldots D_n$ if 
$$\chi_i(P(f,g)) \equiv D_{i,\chi,M}(P)(f,g) \mod h^N$$  
for any $P \in R[[\epsilon]]$ and any positive integer $N$ and $M$
sufficiently large depending on $N$. We also assume that these
equations are true with the convention that $D_0=\partial$ and that
$\chi_0=\chi.$  \endproclaim

\proclaim{Definition 1.3} We say that $D_i$ is slow under the above representation if $D_i(P)$ is not 
in the ideal $(\epsilon) \subset R[[\epsilon]]$ for some $P \in R[[\epsilon]]$, 
but $\chi_i$ vanishes on the 
set $h=0.$\endproclaim

We can construct representations using algebraic geometry. Let $T$ be the disk in $\bold C$ and let $\pi:\Cal X \to T$ be a
smooth proper  family of curves of genus $n.$   We will suppose there are
sections $P:T \to \Cal X$,  $Q:T \to \Cal X$,  $R:T \to \Cal X$ so
that for each $t \in T$, we have $P(t)+Q(t)$ is a divisor on $\Cal X_t
=\pi^{-1}(t)$  linearly equivalent to the divisor $2R(t).$  We will assume that 
$P(0)= Q(0)=R(0)$. Choosing a homology basis of $H_1(\Cal X,\bold Z)$, we can identify $\bold C^g \times T$ with the relative Jacobian of the family $\Cal X \to T.$ Suppose $\gamma \in H_1(\Cal X,\bold Z)$ is given.  We will suppose there is function $h$ on $T$ so that 
$$\int_{Q(t)}^{P(t)} \omega = h(t) \int_{\gamma}\omega\tag 1.3.1$$ for all holomorphic one forms on $\Cal X_t =\pi^{-1}(t).$ Geometrically, (1.3.1) says that under the Abel-Jacobi map which sends a piece of the curve to $\bold C^g$, the secant line from $Q$ to $P$ passes through $\gamma.$ This is meaningful as long as $P(t)$, $Q(t)$ are all close together. Under our identification, a point $(v,t)\in \bold C^g \times T$ gives a line bundle $\Cal L$ on $\Cal X_t.$ When $h(t)=1/N$ for $N$ an integer, we can introduce the functions $A_{\Cal L}$ and $B_{\Cal L}$ as being defined at points $(v,t)$ for $1/t \in \bold Z.$  It is easy to see that there are meromorphic functions $f$ and $g$ which coincide with 
$A$ and $B$, and these $f$, $g$ and $h$ and certain linear flows produce a representation of the Toda $\Cal D_{T_k}.$
Further, by studying the geometry of the situation, we can show that this representation is slow for 
$\Cal D_{T_2}+2\Cal D_{T_1}.$ This turns out to mean that $f$ can be computed asymptotically up to an additive constant from $g$. Lemma 2.6.1. contains the crucial step.  It says that if $f$ is computed in terms of $g$ up to order $h^{n}$ and if the representation is slow, then we can compute $f$ in terms of $g$ up to an additive constant modulo $h^{n+1}.$  This computation is just that there is a $P \in R_0[[\epsilon]]$ so that $f=P(g)+C$ and this $P$ is the desired $P$ of Theorem 1.1. This is a strange relation, since for any particular $N$, there's no relation between $A$ and $B$.  The subtlety here is that the geometric genus of the Bloch spectrum for a arbitrary $A$ and $B$ of periodicity $N$ grows with $N$, but the geometric genus of the Bloch spectrum associated to the $A_{\Cal L}$ and $B_{\Cal L}$ above remains $g$, although the arithmetic genus does grow as we let $t=1/N.$ In fact, $\Cal X_{1/N}$ is the normalization of the Bloch spectrum, but the Bloch spectrum has many nodes which are resolved by the normalization. Further, the line bundle $\Cal L$ on the normalization of the Bloch spectrum becomes a torsion free sheaf on the Bloch spectrum, which is not locally free at new nodes. 

The basic problem turns out to be to construct lots of such representations. We look a families of curves inside $\bold P^1 \times \bold P^1$ of bidegree $(n+1,2)$ with affine coordinates $x$ and $y$. Intuitively, the condition (1.3.1) imposes $g-1$ conditions, so there should be lots of curves satisfying the condition (1.3.1).  We investigate curves near the following curve $C_0$ defined by $$0=(y^2-x)(x-1)(x-\frac 1 {2^2})(x-\frac 1 {3^2})\dots (x-\frac 1 {n^2}).$$ Our object roughly is to show that the subset of curves satisfying (1.3.1) is smooth of codimension $g-1$ near $C_0.$ Further, we let $\Cal L_0$ be a line bundle on $C_0$ of degree $n$ which has degree one on all the vertical components $0=(x-1/k^2)$ of $C_0$ and degree zero on the component $y^2-x=0.$ Then we can explictly calculate the functions $f$ and $g$ we are interested in when we deform the pair $(C_0,\Cal L_0)$ in certain directions. For instance, when we deform one of the nodes of $C_0$ away, but still have the curves satisfying (1.3.1). This gives enough information to produce generic enough $g's.$

There are several technical problems in establishing our results. One is  finding a suitable definition of generic. Another is that $f$ is only determined up to a constant by $g$.  This problem is overcome by a monodromy argument (Lemma 3.7.3).

The sequence of commuting derivations $
\bold D_k$ can be put in the context of Poisson brackets, so that we can consider algebraically the setup of Hamiltonian completely integrable systems with conserved quantities in involution with respect to a Poisson bracket and the associated flows from the conserved quantities. I learned about this type of construction from papers of E. Frenkel.
Let $$\hat R =\bold C[ \ldots \hat a_{-1},\hat a_{0}, \hat a_1 \ldots
\hat b_{-1},\hat b_{0}, \hat b_1 \ldots]$$ We say
a monomial in the $\hat a_k$ and $\hat b_k$ has weight $r$ if the sum of
the subscripts of the  $\hat a_k$ and $\hat b_l$ sum to $r$. So the
monomial $\hat a_1 \hat a_2 \hat b_{-3} $ has weight 0.  Let $I_k
\subset \hat R$ be
the $\bold C$ span of all the elements of weight $k$.  Let $M_N$ be
the ideal of $R$ generated by $$\hat a_N, \hat a_{N+1} \ldots \hat
a_{-N}, \hat a_{-N-1} \ldots
\hat b_N, \hat b_{N+1} \ldots \hat b_{-N}, \hat b_{-N-1} .$$ Let $\hat I_k$ be the completion of $I_k$ with respect to subspaces
$I_k \cap M_N$ as $N \to \infty.$ Then $$\Cal F =\bigoplus_k \hat
I_k$$ is called the Fourier ring.  $\Cal F$ is naturally a graded
ring. We can construct a series of maps $f_n: R[[\epsilon]] \to \Cal F[[\epsilon]]$ so that $f_n(v^{(0)})=\hat a_n$ and $f_n(w^{(0)})=\hat b_n$ and the $f_n$ behave like Fourier coefficients, e.g. $$f_n(HK)=\sum_{k+l=n} f_k(H)f_l(K).$$ One can form an analogous ring $\Cal F_0$ from the ring $$\hat R_0 =\bold C[ \ldots
\hat b_{-1},\hat b_{0}, \hat b_1 \ldots].$$

 One can then show that the Toda derivations on $R[[\epsilon]]$ induce derivations on $\Cal F[[\epsilon]]$ which are compatible with the  $f_n.$  Now each Toda lattice can be put in a Poisson framework and we can make a formal version of these Poisson brackets to obtain a Poisson bracket on $\Cal F[[\epsilon]].$  Further, the Toda flows on  
$\Cal F[[\epsilon]]$  come from conserved quantities in $\Cal F[[\epsilon]].$ Let $\hat I_P\subset \Cal F[[\epsilon]]$
be the closure of the ideal generated by all the Fourier coefficients  $f_n(v^{(0)}-P)$ for $n\neq 0$ and a certain Casimir.  Then we can find an induced Dirac bracket on $ \Cal F[[\epsilon]]/\hat {\Cal I}_P\simeq \Cal F_0[[\epsilon]]$ so that the Toda derivations come by bracketing with conserved quantities. Further, we can find $\hat \beta_k\in \Cal F[[\epsilon]]/(\hat {\Cal I}_P)$ for $k\in \bold Z$ so that modulo $\epsilon$ the $\beta_k$ generate $\Cal F_0$ topologically and satisfy the defining relations of the Virasoro algebra modulo $\epsilon.$

 A similar construction of the deformation of KdV discovered by Frenkel and Reshetikhin \cite{FR} in terms of difference equations has been made by Frenkel \cite{F}. I believe the techniques of this paper will produce  many such deformations of KdV hierarchy as well as deformations of W-algebras.

\head 2. Differential Algebra \endhead

\subhead 2.1 \endsubhead
Let $R$ be the ring of polynomials with complex coefficients with
generators $v^{(i)}$ and $w^{(j)}$ where $i$ and $j$ run over the
non-negative integers,
$$R=\bold  C [v^{(0)},w^{(0)},v^{(1)},w^{(1)},\ldots]$$
We introduce a $\bold  C$ derivation $\partial$ by the formulas 
  $$\partial v^{(i)} = v^{(i+1)}$$ and  $$\partial w^{(i)} = w^{(i+1)}.$$ 
Then $\partial$ on any polynomial in $R$ is defined by the Leibnitz rule. We have a subring $R_0 \subset R$ defined to be the ring generated by the $w^{(n)}.$

This ring $R$ is considered to be the ring of translation invariant
 differential operators in two functions $f(x)$ and $g(x)$.  An
 element of $R$ can be regarded as such a differential operator by
 making the substitutions
$$v^{(n)}= {\partial^n f(x) \over \partial x^n}$$ and $$w^{(n)}=
{\partial^n g(x) \over \partial x^n}$$ so that $\partial$ just becomes
$\partial \over \partial x.$ If $P \in R $ and $f(x)$ and $g(x)$ are
nice functions of $x$, then we define $$P(f,g)(x)$$ to be the result
of making the above substitution.  So for example, if $P =
v^{(1)}w^{(2)}$, then $$P(f,g)(x)= {\partial f(x) \over \partial
x}{\partial^2 g(x)\over \partial x^2}.$$ If $f$ and $g$ depend on a auxiliary variable $t$, then we write $P(f,g)(x,t).$

The ring $R$ can be used to study systems of equations:

$${\partial f(x,t) \over \partial t} = P(f,g)(x,t)$$
$${\partial g(x,t) \over \partial t} = Q(f,g)(x,t),$$
where $P$ and $Q$ are elements of $R$.  We can encode the pair $(P,Q)$ by defining a
derivation $D_{(P,Q)}.$ \definition{Definition 2.1.1} $D_{(P,Q)}$ is the derivation of $R$ commuting with the derivation $\partial$ with the additional properties
$$D_{(P,Q)}(v^{(0)})=P$$ and $$D_{(P,Q)}(w^{(0)})=Q.$$ \enddefinition
Any derivation of $R$ commuting with $\partial$ is of this form.

\subhead 2.2 \endsubhead
We will mostly be concerned with the ring $R[[\epsilon]]$. 
The elements of this ring are formal power series in $\epsilon$ so that the coefficients of 
$\epsilon^n$ are just elements of $R$.  We extend $\partial$ to be a continuous derivation of 
$R[[\epsilon]]$ by taking $\partial \epsilon = 0$.  We next introduce an important series of 
maps $E_k : R[[\epsilon]]\to R[[\epsilon]]$ by the  formulas 
$$E_k(P) = P + k\epsilon \partial P + {k^2\epsilon^2\partial ^2 P \over 2!} + 
{k^3\epsilon^3 \partial ^3 P \over 3!} +\ldots $$
Formally, we can write $$ E_k = \exp(k\epsilon \partial).$$ 
 We have that $$E_k E_j = E_{k+j}.$$
Note that $E_k(v^{(0)})$ is just the Taylor series for $f(x+k\epsilon)$ 
$$E_k(v^{(0)})(f(x),g(x)) = f(x) + {k\epsilon \partial f(x) \over \partial x} + \ldots,$$ when we make the 
substitution of $f(x)$ for $v^{(0)}$ described above.  Note that if $D$
is a continuous derivation of $R[[\epsilon]]$ commuting with
$\partial$ and $D(\epsilon)=0$, then $D$ is uniquely specified by
$D(v^{(0)})$ and $D(w^{(0)}).$ Conversely, given $F$ and $G$ in
$R[[\epsilon]]$, we can find a continuous derivation  $D$ commuting
with $\partial$ and with $D(\epsilon)=0.$ Let us call such a
derivation a tame derivation. 

We will use the ring $R[[\epsilon]]$ to describe the asymptotic behavior of 
difference equations. In our context, a difference equation will be
given by two polynomials $P_1$ and $P_2$ in the ring $S_1$  of polynomials in
the 
variables $\ldots X_{-1},X_0,X_1,X_2 \ldots$ and the variables 
$\ldots Y_{-1},Y_0,Y_1,Y_2 \ldots$. In order to facilitate
substitution, we will write 
$$P_1(\ldots a,\hat b,c \ldots ;\ldots \alpha, \hat \beta
,\gamma,\ldots)$$ to mean the result of substituting $a$ for $X_{-1}$,
$b$ for $X_0$, $c$ for $X_1$ and also substituting $\alpha$ for
$Y_{-1}$, $\beta$ for $Y_0$, etc. That is the $\hat {\phantom{x}}$ is just to
indicate the variable to be substituted for $X_0$ or $Y_0$.  Let $S_2$
be the polynomial ring over $\bold C$ with variables $\ldots
a_{-1},a_0,a_1 \ldots$ and $\ldots
b_{-1},b_0,b_1 \ldots.$ (Of course, $R$, $S_1$ and $S_2$ are all the
same polynomial ring on a denumerable number of variables, but it is
convenient to have different names the variables.)  Given $P_1$ and
$P_2$, we can define a derivation $\Cal D_{P_1,P_2}$ of $S_2$ by
$$\Cal D_{P_1,P_2}(a_n)=P_1( \ldots a_{n-1},\hat a_n,a_{n+1}\ldots;
\ldots b_{n-1},\hat b_n,b_{n+1}\ldots)$$ and $$\Cal D_{P_1,P_2}(b_n)=P_2( \ldots a_{n-1},\hat a_n,a_{n+1}\ldots;
\ldots b_{n-1},\hat b_n,b_{n+1}\ldots).$$ Let $T$ be the automorphism
of $S_2$ defined by $T(a_n)=a_{n+1}$ and  $T(b_n)=b_{n+1}.$ Then $\Cal D_{P_1,P_2}$ is
translation invariant in the sense that $\Cal D_{P_1,P_2}$ commutes
with $T.$  Conversely, any derivation of $S_2$ commuting with $T$ is
of the form $\Cal D_{P_1,P_2}.$ Since the commutator of derivations is
a derivation, this allows us to define the commutator of $(P_1,P_2)$
with $(Q_1,Q_2)$ by 
$$\Cal D_{[(P_1,P_2),(Q_1,Q_2)]}=[\Cal D_{P_1,P_2},\Cal
D_{Q_1,Q_2}].$$

We will now define a continuous derivation $D_{P_1,P_2}$ of
$R[[\epsilon]]$ by 
$$\split D_{P_1,P_2}(v^{(0)})= \\& P_1(\ldots
E_{-1}(v^{(0)}),\hat E_{0}(v^{(0)}),E_{1}(v^{(0)}),\ldots;\\& \ldots,
E_{-1}(w^{(0)}),\hat E_{0}(w^{(0)}),E_{1}(w^{(0)}),\ldots)\endsplit$$
and $$\split D_{P_1,P_2}(w^{(0)})= \\& P_2(\ldots
E_{-1}(v^{(0)}),\hat E_{0}(v^{(0)}),E_{1}(v^{(0)}),\ldots;\\& \ldots,
E_{-1}(w^{(0)}),\hat E_{0}(w^{(0)}),E_{1}(w^{(0)}),\ldots)\endsplit.$$
 We then define  $D_{P_1,P_2}(v^{(n)})=\partial^n D_{P_1,P_2}(v^{(0)})$ and
$D_{P_1,P_2}(w^{(n)})=\partial^n D_{P_1,P_2}(w^{(0)})$ and extend by the
Leibnitz rule. Note that the commutator $[   D_{P_1,P_2},\partial]$
vanishes on the generators $v^{(n)}$ and $w^{(n)}$, so $D_{P_1,P_2}$
commutes with $\partial.$ 
It is an exercise in the chain rule that 
$$[D_{(P_1,P_2)},D_{(Q_1,Q_2)}]=D_{[(P_1,P_2),(Q_1,Q_2)]}.$$

\subhead 2.3 \endsubhead
Suppose that we have an element  $P \in R$.   Let
$\Cal M$ be an analytic manifold and let $\chi$ be an analytic vector
field on $\Cal M$ and let $f$ and $g$ be two meromorphic functions on
$\Cal M$.  We define $P_\chi(f,g)$ to be $P(f,g,\chi f,\chi g, \chi^2
f, \chi^2 g \ldots).$ Next, suppose  we have a function $h$ on $\Cal M.$ For any given $N$, a positive integer, 
we can extend this definition to $P \in
R[[\epsilon]]$ by 
$$P_{\chi,N} (\sum_{n=0}^\infty P_n\epsilon^n)(f,g) =\sum_{n=0}^N
h^nP_{n,\chi}(f,g).$$ If $h_1$ and $h_2$ are two functions on $\Cal M$,
we say 
$$h_1 \equiv h_2 \mod h^N$$ if $(h_1 -h_2)/h^N$ is analytic on an open dense set of the set $\{h=0\}.$ We will use a similar terminology for vector fields on $\Cal M$.

Suppose we are given tame derivations $D_1 \ldots D_n$ of
$R[[\epsilon]]$ and we are given vector fields $\chi,\chi_1 \ldots \chi_n.$
Suppose we are also given a function $h$ on $\Cal M$ which is killed
by $\chi$ and all the $\chi_j.$ We further suppose that the $f$ and
$g$ do not have poles along the set $\{h=0\}$.

\proclaim{Definition 2.3.1} In the above situation, we say that 
$\rho=(f,g,h;\chi,\chi_1 \ldots \chi_n)$  is  a 
representation 
of $D_1 \ldots D_n$ if 
$$\chi_i(P(f,g)) \equiv D_{i,\chi,M}(P)(f,g) \mod h^N$$  
for any $P \in R[[\epsilon]]$ and any positive integer $N$ and $M$
sufficiently large depending on $N$. We also assume that these
equations are true with the convention that $D_0=\partial$ and that
$\chi_0=\chi.$  
We define  $\rho(P)=P(f,g).$  We will only use $\rho(P)$ in congruences modulo $h^N$, so in the context of congruence, the formal power series in $h$ makes sense.
%
Analogously, suppose that $\bar D_1 \ldots \bar D_n \in R_0[[\epsilon]].$ $\rho=(g,h;\chi,\chi_1 \ldots \chi_n)$  is  a 
representation 
of $\bar D_1 \ldots \bar D_n$ if 
$$\chi_i(P(g)) \equiv D_{i,\chi,M}(P)(g) \mod h^N$$  
for any $P \in R_0[[\epsilon]]$ and any positive integer $N$ and $M$
sufficiently large depending on $N$.
\endproclaim
In the Definition, it suffices to check the cases $P=v^{(0)}$ and
$P=w^{(0)}$ to check the equality of the definition for all $P$, 
since both sides are derivations.   
\proclaim{Definition 2.3.2} Suppose $D=\sum a_iD_i$is a linear combinations of the $D_i.$ We say that $D$ is slow under $\rho$ if there is a $P \in R[[\epsilon]]$ so that $D(P)$ is not 
in the ideal $(\epsilon) \subset R[[\epsilon]]$, 
but $\sum a_i\chi_i$ vanishes on the 
set $h=0.$
\endproclaim
 
\subhead 2.4 \endsubhead
We will be constructing representations be in the following context:
Let $V$ be an analytic manifold and let $\Cal M = V \times \bold C^g.$ Let 
$\pi$ be the projection of $\Cal M$ onto $V$.  $h$ will be the pullback of some
function on $V$ via $\pi$.  $\sigma=\sigma_0$ and $\sigma_1\ldots \sigma_k$ will denote 
sections of
$\pi$,  $\sigma_k:V \to \Cal M.$
Now any section $\tau$ of $\pi$ induces vertical vector field
$D_\tau$ on $\Cal M$ by $$\bold D_\tau(f)(x)=\lim_{p\to 0} {f(x+p\tau(x))-f(x)\over
p}.$$

We can define representations of $R[[\epsilon]]$ when we have
$\Cal C^{\infty}$ functions $f$ and $g$ which satisfy difference
equations. \proclaim{Proposition 2.4.1} Suppose that $f$ and $g$ satisfy the
following equations:
$$\split \bold D_{\sigma_n}(f)(x) = \\&P_{1,n}(\ldots
f(x-h(x)\sigma(\pi(x))),
\hat f(x),
f(x+h(x)\sigma(\pi(x))) \ldots;\\&\ldots g(x-h(x)\sigma(\pi(x))),\hat g(x),
g(x+h(x)\sigma(\pi(x))) \ldots)\endsplit $$
$$\split \bold D_{\sigma_n}(g)(x) = \\&P_{2,n}(\ldots
f(x-h(x)\sigma(\pi(x))),
\hat f(x),
f(x+h(x)\sigma(\pi(x))) \ldots;\\&\ldots g(x-h(x)\sigma(\pi(x))),\hat g(x),
g(x+h(x)\sigma(\pi(x))) \ldots)\endsplit. $$
where $P_{l,k}$ are in $S_1$ and the $\hat{\phantom{x}}$ is the place indicator.
We let $D_k$ be the element $D_{P_{1,k},P_{2,k}}$, a tame derivation of
$R[[\epsilon]]$ defined above and let $\chi_k= \bold D_{\sigma_k}$ and
$\chi$ be $\bold D_\sigma$.  Then $(f,g,h;\chi,\chi_1,\ldots \chi_n)$ form a
representation of $D_1,\ldots, D_n.$ \endproclaim
\demo{Proof}

All we need to do is to check that 
$$\chi_i(f) \equiv D_{i,\chi,N}(v^{(0)})(f,g) \mod h^N$$ and 
$$\chi_i(g) \equiv D_{i,\chi,N}(w^{(0)})(f,g) \mod h^N$$

In fact, $$\split D_{i,\chi,N}(v^{(0)})(f,g) \equiv  \\&P_{1,n}(\ldots
f(x-h(x)\sigma(\pi(x))),
\hat f(x),
f(x+h(x)\sigma(\pi(x))) \ldots;\\&\ldots g(x-h(x)\sigma(\pi(x))),\hat g(x),
g(x+h(x)\sigma(\pi(x))) \ldots) \mod h^N\endsplit .$$  
This in turn follows from
$$f(x+p h(x)\sigma(\pi(x))) \equiv E_p(v^{(0)})_{\chi,N}(f,g) \mod
h^N,$$ which in turn is just Taylor's theorem.
\enddemo
\subhead 2.5 \endsubhead
We will frequently use this construction when $$f=-2+h^2 f_1$$ and
$$g=1+h^2 g_1,$$ where $f$ and $g$ are meromorphic functions of $\Cal
M$, which do not have polar divisors containing $\{h=0\}$ To this end,
define a $\bold C$ algebra endomorphism of $R[[\epsilon]]$ commuting
with $\partial$ by $$\Phi(v^{(0)})=-2+\epsilon^2 v^{(0)}$$ and
$$\Phi(w^{(0)})=1+\epsilon^2 w^{(0)}.$$ $\Phi$ does not have an
inverse on $R[[\epsilon]]$, but does have one on $R((\epsilon))$, the ring of
Laurent series in $\epsilon$ which contain only finitely many negative
powers of $\epsilon.$ 

\definition{Definition 2.5.1} Suppose  $D$ is a tame derivation of $R[[\epsilon]]$.
We define a new derivation $D_\Phi$ of $R((\epsilon))$ by
$$D_\Phi = \Phi D \Phi^{-1}.$$ \enddefinition

In situations we will be considering
$D_\Phi$ will turn out to be a tame derivation of $R[[\epsilon]].$ 

\proclaim{Lemma 2.5.2} Suppose $f=-2+h^2 f_1$ and
$g=1+h^2 g_1$. If
$$\rho=(f,g,h;\chi,\chi_1\ldots\chi_n)$$ are a representation of $D_1\ldots D_n$   then 
$\rho_{\Phi}=(f_1,g_1,h;\chi,\chi_1\ldots\chi_n)$ is a representation of 
$D_{1,\Phi}\ldots D_{n,\Phi}$.  \endproclaim
\demo{Proof}  $$(\Phi (P))(f_1,g_1)=P(f,g)$$ for any 
$P\in R[[\epsilon]].$  So 
$$\align \chi_i((\Phi(P))(f_1,g_1))&=\chi_i(P(f,g))
\\&\equiv 
D_i(P)( f,g) \mod h^N\\&\equiv (\Phi D_i(P)(f_1,g_1) \mod h^N
\\&\equiv  (D_{\Phi,i})(\Phi(P))(f_1,g_1) \mod h^N.\endalign$$ 
Given $Q \in R[[\epsilon]]$, we let $P=\Phi^{-1}(Q)$ and then we have 
$$\chi_i(Q)(f_1,g_1))\equiv (D_{\Phi,i})(Q)(f_1,g_1) \mod h^N$$ so we have a
representation. 
\enddemo

Next we work out a simple example of these definitions. We take 
$$P_1=Y_{-1}-Y_0\in S_1$$ and $$P_2= Y_0(X_0-X_1).$$  These are the Toda equations. Then define $D_1$
by  $$\align D_1(v^{(0)})&=
w^{(0)}-w^{(1)}\epsilon+w^{(2)}\epsilon^2/2!+\ldots 
-w^{(0)}\\&=-w^{(1)}\epsilon+w^{(2)}\epsilon^2/2!+\ldots \endalign$$
and $$D_1(w^{(0)})=w^{(0)}(-v^{(1)}-v^{(2)}\epsilon/2!\ldots.$$ 
Then $$D_{1,\chi,N}(v^{(0)})(f,g)=-h\chi( g) +h^2 \chi^2 (g)/2! +\ldots.$$
and   $$D_{1,\chi,N}(w^{(0)})(f,g)=g(-h\chi( f) -h^2 \chi^2 (f)/2!
+\ldots).$$
So if $(f,g,h;\chi_1)$ is a representation of $D_1$, then we will have
$$\chi_1(f)=-h\chi( g) +h^2 \chi^2 (g)/2! +\ldots$$ and 
$$\chi_1(g)=g(-h\chi( f) -h^2 \chi^2 (f)/2!
+\ldots),$$ where these equations are taken to be true  near $\{h=0\}$
modulo high powers of $h$, so the expansions are considered to be
asymptotic. 

\subhead 2.6 \endsubhead
Recall that  $R_0 \subset R$ is the $\bold C$ algebra generated by the
$w^{(n)}$. 
Next, we consider $P \in R_0[[\epsilon]]$ and assume that
$$P=w^{(0)}+\epsilon P_1.$$
Assume we have a representation  $(f,g,h;\chi,\chi_1)$ of $D_1.$
Assume that $$D_1(v^{(0)})=  v^{(1)}-w^{(1)} + \epsilon Q_1 
$$ and
that  $$D_1(w^{(0)}) = w^{(1)}-v^{(1)} +\epsilon Q_2.$$
Assume that this representation  is slow for $D_1$, so that $\chi_1$ vanishes on
${h=0}$. Let $\Cal D$ be the derivation defined by
$$\Cal D(v^{(0)})=  v^{(1)}-w^{(1)} $$ and $$\Cal D(w^{(0)}) = w^{(1)}-v^{(1)}. $$ Then 
$$D_1=\Cal D +\epsilon \Cal E.$$

Let $W =
(v^{(0)}-P)$ and let $\Cal I_P$ be the closure of the ideal generated by $W$, $\partial W$, $\partial^2 W\ldots.$ If $Q\in
R[[\epsilon]]$, let $\bar Q$ be the image of $Q$ in $ R[[\epsilon]]/ \Cal
I_P.$ Then there is the
natural map $$\sigma_1: R_0[[\epsilon]] \to R[[\epsilon]]/ \Cal
I_P.$$ Note that $\sigma_1$ is an isomorphism.  We let $\sigma(Q)$  be
$\sigma_1^{-1}(\bar Q).$  Thus $\sigma(Q)$ is just the result of
replacing any occurrence of $ v^{(0)}$ in $Q$ by $P$, any occurrence of 
$ v^{(1)}$ by $\partial P$, etc.

\proclaim{Lemma 2.6.1}  Given $P$ and $D_1$ as above, there is an $\Cal H\in R_0[[\epsilon]]$ so that if $$ \rho(v^{(0)}-P) \equiv 0 \mod h^n.$$
Then $$ \chi(\rho(v^{(0)}-P))\equiv \rho(\Cal H) \mod h^{n+1}. $$ Note that $\Cal H$ does not depend on $\rho.$\endproclaim

\demo{Proof} 
We let  $s = \rho(v^{(0)}-P) $. Notice that if $Q \in R$, then 
$$ \rho(Q) \equiv \rho(\sigma(Q)) \mod h^n,\tag 2.6.1.1$$ since $P=\sigma(v^{(0)})$
and 
$$ \align\rho(v^{(0)}) &\equiv \rho(P) \mod h^n  \\
&\equiv \rho(\sigma(v^{(0)})) \mod h^n \endalign$$
Let $x=\chi_1(s).$
Then $$\align x&=\rho(D_1(v^{(0)}-P)) \\
&= \rho(D_1(v^{(0)}-w^{(0)} -\epsilon P_1))\\
&= \rho( \Cal D(v^{(0)}-w^{(0)} -\epsilon P_1)+\epsilon \Cal E (v^{(0)}-w^{(0)} -\epsilon P_1))\\
&=\rho(2\partial(v^{(0)}-w^{(0)})+\epsilon( \rho(-\Cal D(P_1)+\Cal E (v^{(0)}-w^{(0)} -\epsilon( P_1))\\
&= 2\chi(\rho(v^{(0)})+\rho(-w^{(0)}))+\epsilon( \rho(-\Cal D(P_1)+\Cal E (v^{(0)}-w^{(0)} -\epsilon P_1))\\
&= 2\chi(s)+2\chi\rho(P)+\rho(-2w^{(1)})+\epsilon( \rho(-\Cal D(P_1)+\Cal E (v^{(0)}-w^{(0)} -\epsilon P_1))\\
&=2\chi(s)+\rho(2\partial P -2w^{(1)})+\epsilon( \rho(-\Cal D(P_1)+\Cal E (v^{(0)}-w^{(0)} -\epsilon P_1))\\
\intertext{Noting that $\rho(\sigma(-\Cal D(P_1)+\Cal E (v^{(0)}-w^{(0)} -\epsilon P_1)))$ is congruent
to $\rho(-\Cal D(P_1)+\Cal E (v^{(0)}-w^{(0)} -\epsilon P_1))$ modulo $h^n$ from (2.6.1.1), we can continue}
&=2\chi(s) + \rho(\sigma(2\partial P
-2w^{(1)}+\epsilon(-\Cal D(P_1)+\Cal E (v^{(0)}-w^{(0)} -\epsilon P_1))\mod h^{n+1}.
\endalign$$
But $$x \equiv 0 \mod h^{n+1}$$ since $D_1$ is slow for $\rho$ so we obtain the conclusion of the Lemma.

\enddemo

\definition{Definition 2.6.2 } 
Let $D_i(v^{(0)})=P_i$ and $D_i(w^{(0)})=Q_i.$
We say $D_1,\ldots D_n$ are  nice if $$\multline  P_i(v^{(0)}+K,v^{(1)},\ldots;w^{(0)},w^{(1)},\ldots)=\\
\sum_{j\leq i} \binom ij
P_j(v^{(0)},v^{(1)},\ldots;w^{(0)},w^{(1)},\ldots)K^j
\endmultline$$
and 

$$\multline  Q_i(v^{(0)}+K,v^{(1)},\ldots;w^{(0)},w^{(1)},\ldots)=\\
\sum_{j\leq i} \binom ij
Q_j(v^{(0)},v^{(1)},\ldots;w^{(0)},w^{(1)},\ldots)K^j
\endmultline
$$ \enddefinition
\definition{Definition 2.6.3}Let $D_i(v^{(0)})=P_i$ and $D_i(w^{(0)})=Q_i.$
 We say the $\rho$ has  weights $r_1\ldots r_n$ if 
$$ P_i(uv^{(0)},uv^{(1)},\ldots;u^2 w^{(0)},u^2 w^{(1)},\ldots)=u^{r_i}P_i(v^{(0)},v^{(1)},\ldots; w^{(0)},w^{(1)},\ldots).$$
and
$$ Q_i(uv^{(0)},uv^{(1)},\ldots;u^2 w^{(0)},u^{2} w^{(1)},\ldots)=u^{r_i+1}Q_i(v^{(0)},v^{(1)},\ldots; w^{(0)},w^{(1)},\ldots).$$

\enddefinition

\definition{Lemma-Definition 2.6.4} Suppose that
$\rho=(f,g,h,\chi;\chi_1\ldots\chi_n)$ is a nice representation and that $s$ is a meromorphic function on $\Cal M$ with 
$\chi(s)=0$ and $\chi_i(s)=0.$  
Let $$\chi_i'= \sum_{j\leq i} \binom ij\chi_j s^j$$ Then 
$\rho'=(f+s,g,h,\chi;\chi_1'\ldots\chi'_n)$ is a representation of  $D_1\ldots
D_n$ and we call $\rho'$  the translation of
$\rho$ by $s.$\enddefinition

\definition{Definition 2.6.5} Suppose $\rho=(f,g,h,\chi;\chi_1\ldots\chi_n)$ has weights $r_1\ldots r_n$.  We define an extended representation $\rho'=(f',g',h',\chi',\chi_1',\ldots)$ on  $\Cal M \times \bold C$ by first defining a function $u:\Cal M \times \bold C\to \bold C$ by $$u(m,z) =1+h(m)^2z$$ and $$f'(m,z)=
u(m,z)f(m)$$ and $$g'(m,z)=u(m,z)^2g(m)$$ 
and $\chi_i'$ is $u^{1-r_i}$ times the natural pullback of $\chi_i$ and $h'$ and $\chi'$ are the natural pullbacks of $h$ and $\chi$ to $\Cal M\times \bold C.$\enddefinition

\remark{Remark 2.6.6} This definition works fine without the particular choice of $u=1+h^2z$ we have made, but the next remark does not.\endremark
\remark{Remark 2.6.7} In the situation of Definition 2.6.5, suppose that $$f=-2+h^2f_1$$ and $$g=1+h^2g_1.$$ Then $$\rho'_{\Phi}(w^{(0)})(m,z)\equiv \rho_{\Phi}(w^{(0)})(m)+2z\mod h.$$\endremark
\definition{Definition 2.6.8} If  $W\in R_0[[\epsilon]],$  we say $g=\rho(
w^{(0)})$ satisfies the equation $W$ non-trivially $\mod h^n$ if 
 $$\rho(W) \equiv 0 \mod h^n$$ and
$g$ does not vanish at a
generic point of $\{h=0\}$ and $W \notin (\epsilon^n).$\enddefinition

\remark{Remark} Note that in the above definition, we can find $W_1 \in
R_0[[\epsilon]]$ so that $W_1 =\epsilon^k W$ and $W_1 \notin
(\epsilon).$ Then $k<n$, so $\rho(W_1)(g)$ is zero when
restricted to  $\{h=0\}.$ \endremark
 \subhead 2.7 \endsubhead
Let $D \in R_0.$  We introduce the variational derivative of $D$
$$ \delta D = \sum_{k=0}^{\infty} (-1)^k \partial^k\left(\frac {\partial
D}{\partial w_k}\right).$$ This operator has the property that 
$$\left[ \frac{d}{d\epsilon}\int_0^{2 \pi} D(f+\epsilon g)
\right]_{\epsilon =0} =\int_0^{2 \pi} \delta   D(f) g,\tag 2.7.0.1$$ if $f$
and $g$ are periodic.
$\delta(D)(f)$ can be somewhat more intuitively defined by equation (2.7.0.1). Evaluating $$\int_0^{2 \pi} D(f+\epsilon g)$$  will yield the integral of a
differential polynomial in $f$ and $g.$ To take the limit at $\epsilon
\to 0$, we can throw away all the non-linear terms in $g$.  Further,
we can eliminate any occurences of derivatives of $g$ by integration
by parts.  The resulting differential polynomial in $f$ is
$\delta(D)(f).$

\proclaim{Lemma 2.7.1} Suppose we have a function $f(t_1,t_2\ldots t_n;x)$
which is periodic with period 1 in $x$.  Suppose that for generic $z$,
we have $$\int_z^{z+1} D(f)(x) dx =0.$$ Then $$\int_z^{z+1}
\delta(D)(f) \frac {\partial f}{\partial t_i} =0,$$ for all
$i$.\endproclaim
\demo{Proof} We assume that $i=1$ for convenience. Fixing $t_1,t_2
\ldots$, we let $$g(x,\epsilon)=\frac { 
f(t_1+\epsilon,t_2 \ldots,t_n)-f(t_1,t_2 \ldots,t_n)}
{\epsilon}$$ Note that $g(x,\epsilon)$ is holomorphic function of $x$,
even when $\epsilon =0.$ In fact, $$g(x,0)=\frac {\partial f}{\partial
 x_1}.$$ So $$\align \int_z^{z+1} D(f+\epsilon g) &= \int_z^{z+1}
D(f(t_1+\epsilon,t_2\ldots t_n;x)dx.  \\
&= 0.\endalign$$ Taking the derivative of both sides of the equation with respect to $\epsilon$
and setting $\epsilon =0$,  we
obtain $$\int_z^{z+1}
\delta(D)(f) \frac {\partial f}{\partial t_1} =0$$

\enddemo
\proclaim{Lemma 2.7.2} Let $D \in R_0$ be of order $m$, i.e. the highest
derivative occuring is of order $m$. Suppose that $$\frac {\partial
f}{\partial t_k}(0,0\ldots 0;x)  =a_k+b_k \exp(2\pi i k  x),$$ with $b_k\neq 0.$ Then
$D(f)$ is not identically zero.\endproclaim
\demo{Proof} Fix a point  $x_0.$ We define a map $\phi_{x_0}:\bold C^m \to \bold C^m$ by
$$\phi_{x_0}\pmatrix  t_1\\t_2 \\t_3 \\ \vdots \\t_m \endpmatrix =
 \pmatrix f(t_1,t_2,\hdots t_m;x_0)\\ \frac{\partial f(t_1,t_2,\hdots
 t_m;x_0)}{\partial x} \\  \frac{\partial^2 f(t_1,t_2,\hdots
 t_m;x_0)}{\partial^2 x}\\ \vdots \\ \frac {\partial^m f(t_1,t_2,\hdots
 t_m;x_0)}{\partial^m x} \endpmatrix.$$ We can compute the Jacobian
 matrix $d\phi_{x_0}$  of $\phi_{x_0}$ at $(0,0,\ldots 0).$
$$d\phi_{x_0}(0,0,0\ldots 0) =
 \pmatrix  a_1+b_1 e^ {2\pi i x_0} & a_2+b_2 e^{ 4\pi i x_0} &\hdots& 
a_m+b_m e^{ 2m\pi i x_0}\\ 2\pi i e^ {2\pi i x_0} & 4\pi ib_2 e^{ 4\pi
 i x_0} &\hdots & 2\pi i m  b_m e^{ 2m\pi i x_0}\\
\hdots & \hdots & \ddots & \hdots \\
 (2\pi i)^m e^{2\pi i x_0} & b_2 (4\pi i)^m e^{ 4\pi
 i x_0} &\hdots & b_m(2\pi i m  )^m e^{ 2m\pi i x_0}
\endpmatrix$$
We claim that $\det(d\phi_{x_0}(0,0,0\ldots 0)) \neq 0$ for generic
values of $x_0.$  It suffices to show that 
$\det(W_{x_0})\neq 0,$ where 
$$ W_{x_0} = \pmatrix \frac {a_1}{b_1}e^{-2\pi i x_0}+1&
 \frac {a_2}{b_2}e^{-4\pi i x_0}+1 &\hdots
 &\frac {a_m}{b_m}e^{-2\pi i mx_0}+1 \\ 
1 & 2 &\hdots &  m \\
1 & 4 & \hdots & m^2 \\
\hdots & \hdots & \ddots & \hdots \\
1 & 2^m &\hdots & m^m
\endpmatrix$$ In particular if $\det (W_{x_0})=0,$  then we would have
$$ \det  \pmatrix \frac {a_1}{b_1}z+1&
 \frac {a_2}{b_2}z^2+1 &\hdots
 &\frac {a_m}{b_m}z^m+1 \\ 
1 & 2 &\hdots &  m \\
1 & 4 & \hdots & m^2 \\
\hdots & \hdots & \ddots & \hdots \\
1 & 2^m &\hdots & m^m
\endpmatrix =0$$ for any value of $z$.  In particular, taking $z=0$,
we would have that a Vandermond determinant was  zero.

So the image of $\phi_{x_0}$ would contain an open set.  On the other
hand, the equation $D(f)=0$ says that  the numbers 
$$\matrix f(t_1,t_2,\hdots t_m;x_0), \frac{\partial f(t_1,t_2,\hdots
 t_m;x_0)}{\partial x} ,  \frac{\partial^2 f(t_1,t_2,\hdots
 t_m;x_0)}{\partial^2 x}, \hdots , \frac {\partial^m f(t_1,t_2,\hdots
 t_m;x_0)}{\partial^m x} \endmatrix $$  satisfy a non-trivial algebraic equation $H=0$
 independent of $t_1, t_2,\hdots t_m.$ But $\{H=0\}$ cannot contain an
 open set.
\enddemo 
\subhead 2.8 \endsubhead
We wish to study representations in the following somewhat degenerate context: We have a
sequence of functions   $g_r$ on $\Cal M_r = V_r \times \bold C$
and let $\pi_r$ be the projection of $\Cal M_r \to  V_r.$
 We can write
$g_r(v,z)$, where $v\in V_r$ and $z \in \bold C.$  
We will assume that  $g_r$ are periodic
with period one 
with respect to the second variable, e.g. $$g_r(v,z) = g_r(v,z+1).$$
Let $\bold D_r$ be the  $\frac {\partial} {\partial z}$
which is tangent to the fibers of $\pi_r.$  We can define $P(g_r)$ for all  $P \in R_0$ so that $(w^{(0)})(g_r)=g_r$ and 
 $$\bold D_r( P(g_r))= \partial(P)(g_r).$$
We also assume there are maps $k_r:V_r \to V_{r+1}$ so that the
pullback  of  $g_{r+1}$ is
$g_r$ via $k_r \times id.$  We also assume we are given
points $Q_r\in V_r$ compatible with the maps $\pi_r.$ 

\definition{Definition 2.8.1} We say the family $\{g_r\}$ is generic if first for
any $D \in R_0$, there is an $r$ so that $D(g_r)$ is not identically
zero.  Second suppose that we are given a  $D \in R_0.$ Each point $v$ of $V_r$ yields a function of $x$ by
$$F_{v,r}(x)=g_r(v,x).$$ Suppose that $$\int_z^{z+1}D(F_{v,r})(x)\,dx=0$$ for generic $z$,$v$ and all $r$. Then there is an $E \in R_0$ so that 
$D = \partial E$.\enddefinition

Next we give a criterion for the family $\{g_r\}$ to be generic. We further
assume that for any positive integer $n$, there is an integer $r$ so that for
any integer $k$, $-n \leq k \leq n$
there are maps $\phi_{k,r}$ of the unit
disk $D \subset \bold C$ to $V_r$ with the following properties: Let
$G_{r,k}$ be the pullback  of $g_r$ to $D \times \bold
C$ via $\phi_{k,r} \times id.$ Let $\bold E$ be the vector field on $D
\times \bold C$
defined by $$(\bold E G)(t,z) = \frac{\partial G(t,z)}{\partial
t}.$$ We suppose that $$(\bold E G_{r,k})(0,z)=a_{k,r}+b_{k,r}\exp(2\pi i k
z),\tag 2.8.1.1$$ where $a_{k,r}, b_{k,r} \in \bold C$ and 
$b_{k,r}\neq 0$ and $a_{r,0}=0.$ Further, we suppose that $\phi_{k,r}(0)=Q_r.$

$\phi_{k,r}$ induces a map $\psi_{k,r}$ from the tangent space $T_0$
of the disk $D$ at zero to the tangent space $T_{Q_r}$ of $V_r.$ Note
that the vectors $$W_{k,r}=\psi_{k,r}(\frac {\partial} {\partial z}),$$ are
all independent, since otherwise we would get a linear dependence
relation between the functions $a_{k,r}+b_{k,r}\exp(2\pi i kz).$ We
can introduce a coordinate system $t_{1},t_{2} \hdots t_p$ on $V_{r}$
so that $Q_r$ is the origin of the coordinate system and the span of
the functions
$$\frac {\partial g_r}{\partial t_k}(0,0,0\hdots 0;x)
$$ contains the span of the functions $\exp(2\pi i kz)$ for $k$ from
$-m$ to $m.$   This is because the linear span of the functions
$\exp(2\pi i k z)$ is the same as the linear span of the functions
$a_{k,r}+b_{k,r}\exp(2\pi i kz).$
\proclaim{Lemma 2.8.2} Suppose that for each $r$, we have that there is a
periodic meromorphic function $w_r$ on $V_r \times \bold C$ so that 
$\bold D_r(w_r)=D(g_r).$  Then under the assumption of (2.8.1.1), there is a $E \in R_0$ so that 
$D = \partial E.$\endproclaim

\demo{Proof} If $\delta(D)=0,$ then we can find $E$ so that $\partial
E = D+C$, where $C$ is constant. But our assumption implies that
$$\int_z^{z+1} D(g_r) =0, $$ for generic $z$, where we integrate along a
straight line for $z$ to $z+1$ in $\bold C.$ Since $E(g_r)$ is periodic,
we have $$\int_{z}^{z+1} \partial E(g_r) =0.$$ So we would have $C=0.$
So we may assume that $\delta(D)\neq 0.$ Note that $\delta(D)\notin \bold C.$ Consequently, we can find a
map $\gamma: D' \to V_m$ so that $$(\delta(D)(g_r))(\gamma(t),x)=t^n
G(t,x),$$ where $G(0,x)$ is not identically zero as a
function of $x$ and $D'$ is the disk. Now $G(v,x+1)=G(v,x).$ Thus for generic $x$ there is an
integer $k$ so that
$$\int_x^{x+1} G(0,z) \exp(2 \pi i k z)dz \neq 0.$$  Choose a $p$ so
that  $\exp(2 \pi i k z)$ is in the span of the functions  
$$\frac {\partial g_p}{\partial t_l}(0,0,0\hdots 0;x)
.$$ On the other hand, we have that 
$$0=\int_z^{z+1} \delta(D)(g_p) \frac {\partial g_p} {\partial
t_l}. $$
Hence we can find a function $R$ so that
 $$0=\int_z^{z+1} \delta(D)(g_p)(\gamma(t),x) R(t,x)\,dx $$ and $$R(Q_p,x)=\exp(2 \pi i k z)$$
 by taking $R$ to be a linear combination of  the functions $
\frac {\partial g_p} {\partial
t_l}. $ In particular, we obtain $$\align 0=&\int_z^{z+1}
 \delta(D)(g_p)(\gamma(t),x))R(t,x) \\
=&t^n \int_z^{z+1}  G(t,x) R(t,x) \endalign$$
So  $$ 0
=\int_z^{z+1}  G(t,x) R(t,x)
.$$  Evaluating at $t=0$ yields a contradiction.
\proclaim{Lemma 2.8.3} Suppose that $D_1 \dots D_n \in R_0.$ and that the family $\{g_r\}$ is generic.  Then we can find an $r$ so that if $D$ is a non-trivial linear combination of the $D_k$ and $D(g_r)=0$, then $D=0.$\endproclaim
\demo{Proof} We can assume that the $D_i$ are linearly independent. Let $W_r \subset \bold C^n$ be the set of $(a_1\dots a_n)$ so that $$\sum_k a_k D_k(g_r)=0.$$ $W_r$ is a sequence of linear subspaces which decrease with $r$.  Further, by the first property of generic $\{g_r\}$, any element of $\bold C^n$ is eventually not in some $W_r.$ \enddemo
 \enddemo
\subhead 2.9 \endsubhead
We will next study the following situation: 
Let $\Cal N_r =V_r\times \bold C^{r}$ where $\bold C^{r}$ has a basis
$\delta_1,\ldots \delta_r$ and let $\pi_r$ be the projection from
$N_r$ to $V_r$.  We suppose there are inclusions $k_r:V_r\to
V_{r+1}.$  For simplicity, we will identify the image of $V_r$ under
$k_r$ with $V_{r}$.  Suppose there are functions $h_r$ on $V_r$ so
that the restriction of $h_{r+1}$ is $h_r$. Let $\Cal Q_r$ be a series of points in $V_r$ with $\Cal Q_r=\Cal Q_{r+1}$ under the identification and assume $h_r(\Cal Q_r)=0$. Then we can send
$\psi_r:\pi^{-1}_{r+1}(V_r) \to N_{r}$ by $$\psi_r(v_r,c_1,\ldots
c_{r+1})=(v_r,c_1,\ldots c_r).$$ Suppose that $f_r$ and $g_r$ are
functions on $N_r$ and that $\psi^*_r(f_r)$ is the restriction of
$f_{r+1}$ and  $\psi^*_r(g_r)$ is the restriction of
$g_{r+1}.$  We will assume that $f_r$ and $g_r$ are generically
defined on the zero section on $\Cal N_r \to V_r.$

Suppose that $D_1 \ldots D_n$ are tame derivations of $R[[\epsilon]]$ and that 
$\chi_{0,r},\chi_{1,r} \ldots \chi_{n,r} $ are vertical vector fields on $\Cal
N_r$  which form representations of $D_1 \ldots D_n$ and so that
$\chi_{0,r}$ represents $\partial.$  We will also assume that $\chi_{0,r}$ is
differentiation in the direction  $\delta_1+2\delta_r+\ldots
r\delta_r.$ Note that $f_{r+1}$ and $g_{r+1}$ are constant
when restricted to the fibers of $\psi_r$ and further we assume that if $F$ is a
function on $\Cal N_r$,
$$\chi_{k,r+1}(\psi^*_r(F))=\psi^*_r(\chi_{k,r}(F))$$ for $k \leq r.$
\definition{Definition 2.9.1} We say that a sequence of representations as
above are a compatible sequence.\enddefinition

Let $W_r\subset \Cal N_r $ be defined by $h_r=0.$ We can generate a
certain class of functions on $W_r$ in the following way: Take any $P
\in R[\epsilon]$ Then $P(f,g)$ is a meromorphic function on $\Cal N_r.$
Then we find the minimal power $p$ of $h_r$ so that $P(f,g)h_r^p$ is generically defined on $W_r.$  We can then restrict $P(f,g)h_r^p$ to $W_r$ to
obtain a meromorphic function on $W_r$.  Let
$\Cal C$ be the set of all functions on
$W_r$ arising from such $P.$ 

\definition{Definition 2.9.2} The sequence of compatible representations is general if the following holds: Suppose  $f$ is in the class $\Cal C.$ Suppose is that $$f_1(v_r,z)=f(v_r,z,2z,\ldots rz)$$ is constant as a
function of $z$. Then $f$ is constant on the fibers of the projection
$W_r \to V_r$. \enddefinition

Now define  functions $f'_r$ and $g'_r$ on $V_r\times \bold C$ by the
formulas:$$f'_r(v_r,z) =f_r(v_r,z,2z,\ldots rz).$$ 
$$g'_r(v_r,z) =g_r(v_r,z,2z,\ldots rz).$$ Notice that the restriction
of $g_{r+1}'$ to $V_r\times \bold C$ is $g_r'.$ Let $\Cal G_{r}$ be the function $g'_r$ restricted to $h_r=0.$ We will assume that the  $(\Cal G_r,\Cal Q_r)$ form a  generic sequence.
We also assume that all functions in the
class $\Cal C$ are
periodic when translated by $\delta_1+2\delta_2+r\delta_r$ and that the compatible sequence of representations is general.

\proclaim{Lemma 2.9.3} Given an integer $n$, 
suppose there is  a $P \in R_0[\epsilon]$ so that \roster
\item $f_r\equiv P(g_r) \mod h_r^n.$ \item $D_1(v^{(0)})\equiv 
v^{(1)}-w^{(1)} \mod \epsilon
$
\item  $D_1(w^{(0)}) \equiv  w^{(1)}-v^{(1)}\mod \epsilon$
\item  $\chi_{1,r} $ vanishes on $\{h_r=0\}$. 
\item $\chi_{0,r} \dots$ form a general representation of $D_1,\dots.$\item The representation is slow for $D_1.$\endroster Then we can find a $Q \in
R_0[\epsilon]$
congruent to $P \mod \epsilon^n$ so that $$\frac {f_r-Q(g_r)}{h^n_r}$$
is constant on the fibers of the projection $W_r \to V_r'.$\endproclaim
\demo{Proof} Using Lemma 2.6.1, we can find a $Q_1 \in R_0[\epsilon]$ so
that $$\chi(f_r) \equiv Q_1(g_r) \mod h_r^{n+1}.$$ In particular, we
can restrict and get $$\bold D_r(f'_r)=Q_1(g'_r). $$ 
Notice that $$(\partial P -Q_1)(g_r') \equiv 0 \mod h_r^n.$$ Since the
$g_r'$ are generic, we must have 
$$\partial P \equiv Q_1 \mod \epsilon^n.$$ If not, we can suppose that
$$\partial P \equiv Q_1 \mod \epsilon^k$$ for some maximal $k < n.$ Then let 
$$E=\frac{ \partial P -Q}{\epsilon^k}$$ and let $F\in R_0$ be the
constant term 
of $E$ as a power series in $\epsilon.$ Then $E(g'_r)$ vanishes on the
set $h_r=0.$ Thus $F=0,$ contradicting the maximality of $k.$ 

So we can write $$Q_1\equiv\partial P
+\epsilon^n Q_2\mod \epsilon^{n+1}$$ for some $Q_2 \in R_0.$   In particular, 
$$h_r^n Q_2( g'_r) \equiv \frac {d}{dz}
\left(f_r'(v_r,z)-P(g_r)(v_r,z)\right) \mod h_r^{n+1} .$$ Thus
$Q_2(
\Cal G_r)$ is the derivative of a periodic function.  Since the
family $\Cal G_r$ is general, we must have $Q_2 = \partial Q_3$ and hence
we may set $Q=P+\epsilon^n Q_3.$ Then $(f_r' -
Q(g_r'))(v_r,z,2z,\ldots rz) $ is constant as a function of $z$.  So
 $$\frac {f_r-Q(g_r)}{h^n}$$ is constant on the fibers of the
projection $W_r \to V_r.$
\enddemo
 
\proclaim{Proposition 2.9.4} Under the hypotheses of Lemma 2.9.3, we can find
functions $F_r$  and vector fields $\chi_1'\ldots \chi'_n$ so
that $F_r,g_r,h_r,\chi,\chi_1'\ldots \chi'_n$ form a representation of
$D_1\ldots D_n$ and a $Q \in R_0[[\epsilon]]$ so that
$$F_r \equiv Q(G_r) \mod h_r^{n+1}$$ and $$Q \equiv P \mod
\epsilon^n.$$ Further,   $\chi_{1,r}' $ vanishes on $\{h_r=0\}$.
\endproclaim
\demo{Proof} Let $s_{n,r}$ be the restriction of $f_r-Q(g_r)$ to the zero
section of $\Cal N_r$ over $V_r$.  We
can think of $s_{n,r}$ as a function on $V_r$ and hence we can think
of  $s_{n,r}$ as a function on $\Cal N_r.$  We define
$F_r=f_r-s_{n,r}.$ We can modify the $\chi_{k,r}$ to $\chi_{k,r}'$ by
translation by $-s_{n,r}$ This new family of representations remains
compatible.   $\chi_{1,r}' $ still vanishes on $\{h_r=0\}$. For 
$$\chi_1 \equiv \chi_1' \mod h_r,$$ since the difference $\chi_1-
\chi_1' $ is divisible by $s_{n,r}$ Note that we start out with $f_r\equiv g_r \mod h_r.$

\enddemo
\proclaim{Theorem 2.9.5}
Suppose there is a generic compatible family of representations with
\roster 
\item $\Cal G_r,\Cal Q_r$ general
\item $\chi_{1,r} $  vanishes on $\{h_r=0\}.$
 \item $D_1(v^{(0)})\equiv 
v^{(1)}-w^{(1)} \mod \epsilon$
\item  $D_1(w^{(0)}) \equiv  w^{(1)}-v^{(1)}\mod \epsilon$\endroster
Then there is a $Q \in R_0[[\epsilon]]$ so that $$D_i(v^{(0)}-Q) \in
\Cal I_Q$$ for all $i.$ 
\endproclaim
\demo{Proof} Suppose we have constructed a $Q_n $ so that $$f_r \equiv
Q_n(g_r) \mod h^n$$ for a given $n$ and all $r$.  By the Proposition 2.9.4, we
can find $$Q_{n+1} \equiv Q_{n} \mod \epsilon^n$$ and modify $f_r$ so that $$f_r \equiv Q_{n+1}(g_r)
\mod h_r^{n+1}.$$ Notice that for any $P \in R[[\epsilon]]$, we have that 
$$P(f_r,g_r) \equiv \sigma_{Q_{n+1}}(P)(g_r) \mod h_r^{n+1}.$$ Since
$\sigma_{Q_{n+1}}$ is a ring homomorphism commuting with $\partial$,
it suffices to check the equation for $P=v^{(0)}$ and $P=w^{(0)}$.
Now $$\align 0 &\equiv \chi_{i,r}(f_r-Q_{n+1}(g_r)) \\&\equiv
D_i(v^{(0)}-Q_{n+1})(f_r,g_r)\\ &\equiv
\sigma_{Q_{n+1}}(D_i(v^{(0)}-Q_{n+1}))(g_r)
\mod h_r^{n+1}.\endalign$$  We let $Q$ be the limit of the $Q_n.$ Then
$$\sigma_{Q}(D_i(v^{(0)}-Q))(g_r) \equiv 0 \mod h_r^{n}$$ for
all $n$.  Since the $g_r$ are generic, we have
$$\sigma_{Q}(D_i(v^{(0)}-Q)) =0, $$
i.e. $D_i(v^{(0)}-Q) \in \Cal I_Q.$
 
\enddemo
\definition{Definition 2.9.6} Suppose that $$D_i(v^{(0)}-Q) \in
\Cal I_Q.$$  Then we define  derivations $$\bar D_i:R_0[[\epsilon]]\to R_0[[\epsilon]]$$ by $$\bar D_i(P)=\sigma(\overline {D_i(P)})$$ where 
$\sigma$ is the inverse of the natural map $$\sigma_1:R_0[[\epsilon]]\to R[[\epsilon]]/\Cal I_Q$$ and $\bar T$ is the image of $T\in R[[\epsilon]]$ in $R[[\epsilon]]/\Cal I_Q.$
\enddefinition
\proclaim{Theorem 2.9.7} Suppose that $\rho=(f,g,h;\chi_1,\ldots \chi_n)$ is a representation of
$D_1,\ldots D_n$ and that $$f \equiv Q(g) \mod h^n$$ for all $n.$  Then using the notation of Definition 2.9.6, $(g,h;\chi_1\dots \chi_n$ is a representation of $\bar D_1 \dots \bar D_n.$\endproclaim
\demo{Proof} 
$$\align \chi_i(g) &\equiv D_i(w^{(0)}(f,g)) \mod h^n \\
&\equiv D_i(w^{(0)}(Q(g),g)\mod h^n\\
&\equiv \sigma(D_i(w^{(0)})(g)\mod h^n\\
&\equiv \bar D_i(w^{(0)})(g)\mod h^n
\endalign   $$ \enddemo
\definition{Definition 2.9.8} If $D \in R_0[[\epsilon]]$ so that $$D=D_0+D_1\epsilon+D_2\epsilon^2+\ldots,$$ we define
$$D^{[0]}=D_0.$$ Similarly, we let $\chi_i^{[0]}$ be the restriction of  $\chi_i$ to the set $\{h=0\}$.\enddefinition
\proclaim{Lemma 2.9.9} Suppose that $D_1 \ldots D_n \in R_0.$ Suppose that $E_1\ldots E_n \in R_0.$ Suppose that $(g,0;\chi_1,\ldots \chi_n)$ is a representation of $D_1,\ldots D_n$ and $(g,0;\chi'_1,\dots \chi'_n)$ is a representation of $E_1\dots E_n.$ Assume further that $\chi_1\dots \chi_n$ and $\chi'_1\dots \chi_n'$ span the same $n$ dimensional vector space.  Further, assume that if $g$  satisfies any non-trivial linear combination $D$ of the $D_i$ and  $E_i$, then $D=0$. Then  linear span of the $D_i$ is the same as the linear span of the $E_i.$ \endproclaim
\demo{Proof} We can assume that the $\chi'_i=\chi_i.$ Then $(D_i-E_i)(g)=0$, so $D_i=E_i.$\enddemo
\proclaim{Corollary 2.9.10} Suppose that $D_1 \in R_0.$ Suppose that $(g,0;\chi_1)$ is a representation of $D_1.$ Further, assume that $D_1(g)=0$ implies that $D_1=0.$  Then if $\chi_1(g)=0$, then $D_1=0.$\endproclaim

\head 3. Constructing  Deformations\endhead

\subhead 3.1 \endsubhead
Let $T$ be a smooth analytic manifold and let $\pi:\Cal X \to T$ be a
proper flat family of curves of arithmetic genus $n.$  Let $t_0 \in T$
be a fixed point. We will usually be interested in the behavior of the
family and related objects near $t_0$. We will suppose there are
sections $P:T \to \Cal X$,  $Q:T \to \Cal X$,  $R:T \to \Cal X$ so
that for each $t \in T$, we have $P(t)+Q(t)$ is a divisor on $\Cal X_t
=\pi^{-1}(t)$ is linearly equivalent to the divisor $2R(t).$ We will
also assume that $P(t)$, $Q(t)$ and $R(t)$ are all smooth points of
$\Cal X_t$.  

Let $\delta:T\times S^1 \to \Cal X$ be a map defined over $T$ so that
$\pi$ is smooth at 
every point of the image of $\delta$ and let $\omega$ be a section of
$\pi_*(\omega_{\Cal X/T}).$ ($\omega_{\Cal X/T}$ is the sheaf of relative
dualizing differentials.)  Then $\delta(t)$ is a
cycle on $\Cal X_t$ and we can form 
$$I_{\delta}(t)(\omega_t)=\int_{\delta(t)} \omega_t$$ where $I_{\delta}$ is
a section of the dual of $\pi_*(\omega_{\Cal X/T})$.  So we can
consider  $I_{\delta}$ as a section of $R^1\pi_*(\Cal O_{\Cal X}).$ We
will suppose that we can find $\delta_1 \ldots \delta_n$ so that
$I_{\delta_1} \ldots I_{\delta_n}$ are a basis of  $R^1\pi_*(\Cal
O_{\Cal X}).$ Let $v_k$ be sections of $\pi_*(\omega_{\Cal X/T})$
which are dual to the $I_{\delta_k}.$ 

We will assume that there is an open set $U \subset \Cal X$ so that $z:U
\to D \times T$ is an isomorphism of  $D \times T$ as $T$ spaces, where $D$ is the
unit disk.  We assume that  the images of the sections $P$, $Q$ and
$R$ are all contained in $U.$ For simplicity, we will assume that the
image of $P$ is the zero section and let $$h =z(R).$$ We will assume that $$dz=v_1.$$ Further, we will assume that we can find anti-derivatives $h_k$ of the $v_k$ defined on $U$ and that the functions 
$$\frac{\int_Q^P v_k}{\int_Q^P v_1}=k.$$ This condition insures  that the secant line between $Q$ and $P$ on the local curve $$t \to (h_1(t),h_2(t),\ldots h_n(t))$$ passes through $(1,2,\ldots,n)$ when $h_1(t)\neq 0$ and insures that the tangent line to $P$ passes through $(1,2,\ldots,n)$ when $h_1(t)=0.$  The existence of such families is not obvious at this point, but we will exhibit such families in \S 4.

\subhead 3.2 \endsubhead
Let $$T_1 = T \times \bold C^n$$ and let $$\Cal X_1 = \Cal X \times
\bold C^n =  \Cal X \times_T T_1.$$  We will define a
relative line bundle $\Cal L$ on $\Cal X_1$ over $T_1$. Let $\pi_1:\Cal X_1 \to T_1$ be the projection. We have the coordinate
functions $z_1 \ldots z_n$ on $\bold C^n.$ We will again denote the pullback
of $z_i$ to $T_1$ by  $z_i.$  On the other hand, we will denote the
pullback of $I_{\delta_i}$ to a section of $R^1\pi_{1,*}(\Cal O_{\Cal X_1})$
by   $I_{\delta_i}$.  Thus we can consider $I=\sum z_k I_{\delta_k}$
as a section of 
$R^1\pi_{1,*}(\Cal O_{\Cal X_1})$. Then there is a relative line bundle
$\Cal L$ on $\Cal X_1$ so that $\Cal L$ corresponds to the same section of
$R^1\pi_{1,*}(\Cal O_{\Cal X_1}^*)$ as  $\exp(2\pi i I).$  We also assume that
we have a line bundle $\Cal M$ of relative degree $n$ on $\Cal X_1$ which is the pullback of a line bundle on $\Cal X \to T.$ Let
$\Cal N = \Cal M \otimes \Cal L.$

We will assume that $P(T)$ and $Q(T)$, which are divisors on $\Cal X,$ meet
transversally. Consider  $Z =\pi(P(T)\cap Q(T)),$ where $\pi$ is the projection of $\Cal X\to T.$  
We will assume $Z$ is a divisor in $T$ and that the map from $ P(T)\cap Q(T)\to Z$ is an isomorphism. Let $Z_1$ be the inverse image of $Z$ in $T_1.$ We
assume that $R(T)\cap P(T)$ and $R(T)\cap Q(T) $ are both equal to $P(T)\cap Q(T). $ We will assume that $Z_1$ is defined by an
equation $h=0$, where $h$ is the pullback of a function on $T.$ 

We will also assume that we have chosen a 
non-constant
section $\lambda$ of $\Cal O_{\Cal X}(P(T)+Q(T)-2R(T))$   Let's fix a
point $t\in T_1$.  Then on $\Cal X_{1,t}$, we have the line bundles
$\Cal N_{k,t}=\Cal N_t(k(P(t)-Q(t)))$.
\definition{Definition 3.1.1} Suppose $P(t)\neq Q(t).$ Then $t$ is $N$-good if each of the line bundles $\Cal
N_{k,t}(-P(t))$, is
non-special for $|k|<N$. We then define sections $s_k\in H^0(\Cal N_{k,t})$ for $|k|<N$ so that
$$\frac {s_{k+1}}{\lambda s_k}(P(t))=1.$$\enddefinition
\remark{Remark} Given our choice of $\lambda$, these $s_k$ are determined up to an
non-zero multiplicative constant independent of $k.$\endremark
We will assume we have made a special choice to $\lambda$ to be denoted $\lambda_0$ so that $z\lambda_0(P)=1.$ Let's fix $N$ and let $T_2 \subset T_1$ be the set of $N$-good points.  
\proclaim{Theorem 3.1.2}  
There are functions $A_{k,t,\lambda}$ and 
 $B_{k,t,\lambda}$ defined on $T_2$ for $|k|<N-1$ so that $$\lambda s_{k}=s_{k+1}+A_{k,t,\lambda}
s_{k}
+B_{k,t,\lambda} s_{k-1}.$$ The function $B$ is never zero on $T$. \endproclaim

\demo{Proof} The dimension of $H^0(\Cal
N_p(P(p)+Q(p)))=3.$  So there must be a linear dependence
relations $$C\lambda s_{k}+Ds_{k+1}+As_{k}+Bs_{k-1}=0,$$ where
$A$, $B$, $C$ and $D$ are all in $\bold C.$ Now we have chosen our
normalizations of the $s_m$ and of $\lambda$ so that 
$$ \frac {s_{k+1}}{\lambda s_{k}}=1$$ at $P(p)$. We must then have 
$C=-D.$  Further, we cannot have $C=D=0$, since this would lead to a
dependence relation between  $s_{k}$ and $s_{k-1},$
which have different order poles at $P(p)$. So we can completely
normalize by taking $C=-1$ and $D=1$ and finally writing 
$$\lambda s_{k}=s_{k+1}+A s_{k}+B s_{k-1}.$$
\enddemo 

If $u$ is a nowhere zero function on $T$, then 
$$A_{k,t,u\lambda}=uA_{k,t,\lambda}$$ and 
$$B_{k,t,u\lambda}=u^2A_{k,t,\lambda}.$$  Further, we have 
$$A_{k,t,C+\lambda}=A_{k,t,\lambda}+C,$$ while $B$ is unchanged by
adding $C$ for any function $C$ on $T.$ We can now normalize the $B_k$ in the following way: Let $$\bold B_{k,t} = \frac{B_{k,t,\lambda}}
{B_{k,0,\lambda}}.$$

\subhead 3.3 \endsubhead
We next discuss theta functions following [M]. Let $C =\Cal X_t$ and denote $P_t$ by $P$, etc. We
assume that $C$ is non-singular. We regard $H^1(C,\bold
Z) $ as a subgroup of $H^0(C,\Omega)^*$. Let $\bold C_1^*$ be the set
of all complex numbers of absolute value one.  Choose a map $$\alpha:
H^1(C,\bold Z) \rightarrow \bold C_1^*$$  so that $$ {\alpha(u_1+u_2) \over
\alpha(u_1)\alpha(u_2)}=e^{i\pi \langle u_1,u_2 \rangle}.$$  There is
a unique hermitian form $H$ on $H^0(C,\Omega)^*$ so that $${ \Im}\,H(x,y)=
\langle x,y \rangle.$$  
Let $\vartheta$ defined on $H^0(C,\Omega)^*$ be the function
satisfying the functional equation
$$\vartheta(z+u)=\alpha(u)e^{\pi H(z,u)+\pi H(u,u)/2}\vartheta(z)$$ for $z \in
H^0(C,\Omega)^*$ and $u \in H^1(C,\bold Z)$. There is a map
$\gamma:U\cap C\to H^0(C,\Omega)^* $ defined by $$\gamma(q) =\int_{P}^q
\omega, $$ where the path from $P$ to $q$ is chosen to lie in $U\cap C$.
Thus we can regard $\gamma(q)\in H^1(C,\Cal O).$ Given any
non-special line bundle $\Cal L=\Cal O_{C}(D)$ of degree $n$ on $C$
with $D$ effective, here is a constant $K_{\Cal L}\in H^1(\Cal O_C)$
so that the zeros of the function  $ p\to \vartheta(\gamma(p)+K_{\Cal L})$ are just $D$ counting multiplicity. Further, $$\exp(2\pi i
K_{\Cal L})=\Cal L$$ Also $$K_{\Cal L(q-P)}=K_{\Cal L}-\gamma(q)$$
modulo periods. Let $K_0$ be a constant corresponding to the line
bundle $\Cal M$ and choose a line bundle $\Cal M_1=\Cal O(D_1)$ so
that $P$ is a  point of multiplicity one of $D_1$ and all the other
points of $D_1$ are outside $U$. Select $K_1$
corresponding to $\Cal M_1.$ So $\vartheta(\gamma(z)+K_1)$ vanishes
exactly once at $P$ and at no other point of $U.$

Using the theta function, we can write down an expression for
a function $\lambda_0$ initially valid in $U$. 
$$\lambda_0(z)=\alpha \frac{\vartheta(\gamma(z)-
\gamma(R)+K_1)^2}{\vartheta(\gamma(z)+K_1)\vartheta(\gamma(z)+K_1-\gamma(Q))}
,$$ where $$\alpha = \frac{\nabla(\vartheta)(K_1)\cdot \gamma'(P)\vartheta(K_1-\gamma(Q))}{\vartheta(K_1-\gamma(R))^2}.$$
This is a well-defined meromorphic function on $C$, since $P+Q$ is
linearly equivalent to $2R$ so by Abel's theorem,
$\gamma(P)+\gamma(Q)=2\gamma(R).$ Using the periodicity properties of
$\vartheta$,  $$\frac{\vartheta(Z-
\gamma(R)+K_1)^2}{\vartheta(Z+K_1)\vartheta(Z+K_1-\gamma(Q))}$$ is
periodic in $Z \in H^1(\Cal O).$

Next, we develop a formula for $$u=\frac {\lambda_0 s_0}{s_1}.$$  Attached to the
point $t$, there is a line bundle $\Cal N_{0,t}$ and let $L$ be the
projection on $t$ to $\bold C^n.$  The zeros of
$\vartheta(\gamma(z)+L+\gamma(Q))$ match the zeros of $s_1$ and  the zeros
of $\vartheta(\gamma(z)+L)$ match the zeros of $s_0$. On the other
hand,
 $$\frac{\vartheta(\gamma(z)+K_1-\gamma(Q))}{\vartheta(\gamma(z)+K_1)}$$ has a simple pole at $P$ and a simple zero at $Q.$  So $s_1/s_0$ is a multiple of the
rational function
$$  \frac{\vartheta(\gamma(z)+K_1-\gamma(Q))
\vartheta(\gamma(z)+L+\gamma(Q))}
{\vartheta(\gamma(z)+K_1)\vartheta(\gamma(z)+L)}.$$
Consequently, $u$ is a multiple of 
$$ \frac{\vartheta(\gamma(z)+K_1-\gamma(R))^2\vartheta(\gamma(z)+L)}
{\vartheta(\gamma(z)+K_1-\gamma(Q))^2\vartheta(\gamma(z)+L+\gamma(Q))}.$$
Using that $u(P)=1$ and $\gamma(P)=0$ we obtain that 
$$u= \frac
{\vartheta(\gamma(z)+K_1+\gamma(R))^2\vartheta(\gamma(z)+L)\vartheta(K_1-\gamma(Q))^2
\vartheta(L+\gamma(Q))}{\vartheta(\gamma(z)+K_1-\gamma(Q))^2
\vartheta(\gamma(z)+L+\gamma(Q))\vartheta(K_1-
\gamma(R))^2\vartheta(L)}.$$ 

Let $a$ denote the constant term of the Laurent series for $\lambda -
s_1/s_0$ developed around $z=0.$ We get $a=\frac{d}{dz}(u)$ evaluated at $z=P.$
Since $u=1$ at $z=0$, we obtain 
$$\split a=-2\frac{\nabla(\vartheta)(K_1-\gamma(Q))\cdot \gamma'(P)}
{\vartheta(K_1-\gamma(Q))}-
\frac{\nabla(\vartheta)(L+\gamma(Q))\cdot \gamma'(P)}
{\vartheta(L+\gamma(Q))}+\\ 2\frac{\nabla(\vartheta)
(K_1+\gamma(R))\cdot \gamma'(P)}
{\vartheta(K_1+\gamma(R))}+\frac{\nabla(\vartheta)(L)\cdot \gamma'(P)}
{\vartheta(L)}\endsplit \tag 3.1.2.1$$

 For future reference, we give  a formula for the constant term $C_0$
of the Laurent expansion of $s_1/s_0.$, namely 
$$\split C_0= \frac{\nabla(\vartheta)(K_1-\gamma(Q))\cdot \gamma'(P)}
{ \vartheta(K_1-\gamma(Q))}+ 
\frac{\nabla(\vartheta)(L+\gamma(Q))\cdot \gamma'(P)}
{\vartheta(L+\gamma(Q))}-\\ \frac{\nabla(\vartheta)(L)\cdot \gamma'(P)}
{\vartheta(L)}-\frac{W}{\nabla(\vartheta)(L)\cdot
\gamma'(P)},\endsplit \tag 3.1.2.2$$ where
$$W=\frac 1 2\frac{d^2}{dz^2}(\vartheta(\gamma(z)))$$ evaluated at $z=P.$

We will also be interested in evaluating $$v=\frac{\lambda_0
s_{0}}{s_{-1}}$$ at $Q$. We obtain that  
$$\frac{s_{-1}}{s_{0}}=
F\frac{\vartheta(\gamma(z)+K_1)\vartheta(\gamma(z)+L-\gamma(Q))}
{\vartheta(\gamma(z)+K_1-\gamma(Q))\vartheta(\gamma(z)+L)}$$
with $$F=\frac{\vartheta(K_1-\gamma(Q))\vartheta(L)}
{\nabla(\vartheta)(K_1)\cdot
\gamma'(P)\vartheta(L-\gamma(Q))}$$ So
$$\frac{\lambda_0 s_{0}}{s_{-1}}= F_1
\frac{\vartheta(\gamma(z)+L)\vartheta(\gamma(z)-\gamma(R)+K_1)^2}
{\vartheta(\gamma(z)+K_1)^2\vartheta(\gamma(z)+L-\gamma(Q))}$$, where 
$$F_1=\frac{(\nabla(\vartheta)(K_1)\cdot \gamma'(P))^2\vartheta(L-\gamma(Q))}
{\vartheta(-\gamma(R)+K_1)^2\vartheta(L)}.$$
So evaluating at $Q$, we get 
$$\align b=&\frac{\lambda_0 s_{0}}{s_{-1}}(Q)\\=& F_1
\frac{\vartheta(L+\gamma(Q))\vartheta(\gamma(Q)-\gamma(R)+K_1)^2}
{\vartheta(\gamma(Q)+K_1)^2\vartheta(L)}\\
=&\left(\frac{\vartheta(\gamma(Q)-\gamma(R)+K_1)\nabla(\vartheta)(K_1)\cdot
\gamma'(P)}
{\vartheta(-\gamma(R)+K_1)\vartheta(\gamma(Q)+K_1)}\right)^2
\frac{\vartheta(\gamma(Q)+L)\vartheta(-\gamma(Q)+L)}{\vartheta(L)^2}\\
=& \left(\frac{\vartheta(\gamma(R)+K_1)\nabla(\vartheta)(K_1)\cdot
\gamma'(P)}
{\vartheta(-\gamma(R)+K_1)\vartheta(2\gamma(R)+K_1)}\right)^2
\frac{\vartheta(2\gamma(R)+L)\vartheta(-2\gamma(R)+L)}{\vartheta(L)^2}\tag 3.1.2.3 \endalign$$
bearing in mind that $2\gamma(R)=\gamma(Q).$

We have defined $a$ and $b$ in equations (3.1.2.1) and (3.1.2.3). $a$ and $b$ depend on $t.$

\proclaim{Proposition 3.1.3} For our choice of $\lambda_0$, we have 
$$a=A_{0,t,\lambda_0}$$ and 
$$b=B_{0,t,\lambda_0}$$ near $t_0.$\endproclaim

\subhead 3.2 \endsubhead
We
will now make a canonical choice of $\lambda$ near the set $Z.$
Specifically, let $$U = 
\sqrt{ B_{0,(x,0),\lambda_0}}.$$ The expression for $b$ makes it clear
that such a meromorphic function exists, since $$\lim_{\gamma(R) \to 0}
\frac{\vartheta(L)^2}{\vartheta(2\gamma(R)+L)\vartheta(-2\gamma(R)+L)}=1$$
for any $L$, in particular for $L=0.$  In particular, we set 
$$U=\left(\frac{\vartheta(\gamma(R)+K_1)\nabla(\vartheta)(K_1)\cdot
\gamma'(P)}
{\vartheta(-\gamma(R)+K_1)\vartheta(2\gamma(R)+K_1)}\right)\sqrt{\frac{\vartheta(2\gamma(R))\vartheta(-2\gamma(R))}{\vartheta(0)^2}}$$
where we choose the branch of the square root which is near 1. When $R$ is near $P$, then $$\vartheta(\gamma(R)+K_1) \approx \gamma(R)\cdot \nabla(\vartheta)(K_1)$$ and $$\vartheta(-\gamma(R)+K_1) \approx -\gamma(R)\cdot \nabla(\vartheta)(K_1)$$ and 
$$\vartheta(2\gamma(R)+K_1) \approx 2\gamma(R)\cdot \nabla(\vartheta)(K_1)$$ so $$U \approx -\frac{\nabla(\vartheta(K_1)\cdot \gamma'(P)}{
\nabla(\vartheta)(K_1)\cdot \gamma(Q)}.$$ On the other hand, we have that $$\vartheta (K_1 -\gamma(Q))\approx -\nabla(\vartheta(K_1))\cdot \gamma(Q),$$ so the
constant term of $s_1/s_0$ is approximately $$
 -\frac{\nabla(\vartheta)(K_1)\cdot \gamma'(P)}{
\nabla(\vartheta)(K_1)\cdot \gamma(Q)}$$ from the expression for $C_0.$ So the constant term of
 $U^{-1} s_{1}/s_0$ is just $1$. Now choose $C$, a function on $T$,
 so that $$U^{-1}A_{0,(x,0),\lambda_0}+C=-2.$$

\proclaim{Proposition 3.2.1} We can canonically choose $$\lambda_{can}=U^{-1}\lambda_0+C$$ so that \roster
\item $2+A_{k,t,\lambda_{can}}$ is divisible by $h^2$.
\item $-1+B_{k,t,\lambda_{can}}$ is divisible by $h^2$.
\item  $2+A_{k,0,\lambda_{can}}=0$ 
\item  $-1+B_{k,0,\lambda_{can}}=0$ 
\item  We have $$\lambda = \alpha_{-1}\frac h z +\alpha_0+\alpha_1\frac z h +\ldots.$$  where $\alpha_{-1}(t)$ vanishes on $Z$ and $\alpha_0$ becomes $-1$  on
 $Z$.\endroster
Of course, the canonical choice depends on the line bundle $\Cal M,$
 but is otherwise completely determined near $Z$ by properties three,
 four and five.   
 \endproclaim
\demo{Proof} Note that $$V(z,h)=\lambda_0 \frac {z(z-Q)}{(z-R)^2}$$ does not have zeros or poles when $h\neq 0$ and $h$ and $z$ are  near zero. 
 Now near $h=z=0$, $V(z,h)$ is a power of $h$ times a unit. We have $$\lambda_{can}=\frac {U^{-1}s_{n+1}}{s_n}+A_{can}+B_{can}\frac{Us_{n-1}}{s_{n}}$$ Hence the constant term of $\lambda_{can}$ as a Laurent series in $z$ is the constant term of $$\frac {U^{-1}s_{n+1}}{s_n}$$ plus the constant term of $A_{can}$.  So modulo $h$, the constant term of $\lambda_{can}$ is $-1$. The Laurent series for $\lambda_0$ is $$\frac 1 z +\frac {J}{h}+\ldots.$$ So $hV(c,h)$ is a unit. But the Laurent series of $hV(z,h)$ is a power series in $z/h.$ So $\lambda_{can}$ is a power series in $\frac z h.$\enddemo
\definition{Definition 3.2.2} $$\bold A_{k,t}=A_{k,t,\lambda_{can}}$$
$$\bold f = \bold A_{0,t}$$
  $$\bold g = \bold B_{0,t}$$ $\bold f$ and $\bold g$ are functions on $T\times \bold C^n.$
\enddefinition

\proclaim{Lemma 3.2.3} $\bold f$ and $\bold g$ are  defined near any 2-good
point.\endproclaim
Let $\chi$ be differentiation in the direction $(1,2,\ldots n).$ 
\proclaim{Proposition 3.2.4} $$\bold g(t,L)=1+4h(t)^2(\chi^2(\log(\vartheta)(L)-\chi^2(\log(\vartheta)(0))+\text{higher order terms in h}.$$\endproclaim
\demo{Proof} We define $$\bold g_1(L) = \frac{\vartheta(L)^2}{\vartheta(L+\gamma(Q))(\vartheta(L-\gamma(Q))}$$ and then $$\bold g =\frac{\bold g_1(L)}{\bold g_1(0)}.$$ On the other hand, $$\gamma(Q)=2h(1,2,\ldots,n)$$ and we know that $\gamma''(Q/2)=0.$ \enddemo
Let's examine the function
$\bold G(t,L)$ which is the analytic continuation of  $$\frac {\bold g(t,L)-1}{h(t)^2}.$$  We can introduce a series of vectors $\bold v_k \in \bold C^n$ by using the identification of $\bold C^n$ with $$H^0(\omega_{\Cal X_{1,t}})^*$$ via integration over the $\delta_i$ by $$\bold v_k(\omega) =Res_{P(t)} \frac{\omega}{z^{k}}.$$ Let $\bold K_1,\bold K_2 \ldots$ be the KdV differential operators. Let $\chi_k$ be directional derivative in the directions $\bold v_k.$
\proclaim{Proposition 3.2.5} We have 
$$\lambda = \alpha_{-1}\frac h z +\alpha_0+\alpha_1\frac z h +\ldots.$$ We can find $\beta_{k,l}$ which are universal polynomials in the $\alpha_k$ so that $$\sum_{l\leq k} \beta_{k,l}Res_{P(t)}\lambda^l \omega \equiv  h^{2k-1}\bold v_{2k+1}(\omega) \mod h^{2k}.$$ Further, $$\beta_{k,k}=n_k\alpha_{-1}^{n_k},$$ where $n_k\in \bold Q $ and $n_l\neq 0$\endproclaim
\demo{Proof} First note that we can find $\gamma_{k,l}$, universal polynomials in the $\alpha_k$ so that
$$\sum_{l\leq k} \gamma_{k,l}Res_{P(t)}\lambda^l \omega \equiv \alpha_{-1}^k \frac {h^k}{z^k}.$$ Indeed, $$\lambda = \alpha_{-1}\frac h z +\ldots.$$ $$\lambda^2=\alpha_{-1}^2\frac {h^2}{z^2} +2\alpha_{-1}\alpha_{0}\frac h z+\ldots.$$
$$\lambda^3=\alpha_{-1}^3 \frac {h^3}{z^3}+3\alpha_0\alpha_{-1}^2\frac {h^2}{z^2}+(3\alpha_{-1}\alpha_{0}^2+3\alpha_{-1}^2\alpha_1)\frac h z+\ldots$$ So we can  express $$\mu_k=\alpha_{-1}^k \frac {h^k}{z^k}$$ in terms of $\lambda \ldots \lambda^k$ with the coefficients universal polynomials in the $\alpha_k$'s.

Note that $z-h$ vanishes on $R(t)$ and that $z-h$ is the anti-derivative of a holomorphic differential.  Since the curve is hyperelliptic and $R(t)$ is a Weirstrass point, $z-h$ is an odd function under the involution of the hyperelliptic curve.  However, any differential is even under the involution.  Thus any differential can be expanded around $P(t)$ as a power series in the even powers of $z-h$ times $dz.$  So if $\omega$ is a differential,  we can write $$ \align \omega =& (a_0+a_2(z-h)^2+a_4(z-h)^4 \ldots)dz\\=&f((z-h)^2)dz.\endalign$$ So $$Res_{P(t)}\mu_1\omega = \alpha_{-1}hf(h^2).$$ $$Res_{P(t)}\mu_2\omega =2h^3\alpha_{-1}^2f'(h^2).$$
$$Res_{P(t)}\mu_3\omega =h^3\alpha_{-1}^3(2f'(h^2)+4h^2 f''(h^2)).$$
$$Res_{P(t)}\mu_4\omega =h^4\alpha_{-1}^4(2hf''(h^2)+8hf''(h^2)+8h^3f'''(h^2)).$$ Continuing in this way, we see that 
we can express $$h^{2k-1}f^{(k)}(h^2)$$ as a linear combination of $$Res_{P(t)}\mu_l\omega $$ for $l\leq k.$ On the other hand,$$f^{(k)}(h^2)\equiv \bold v_{2k+1}(\omega)\mod h.$$

\enddemo
\definition{Definition 3.2.6} Suppose $D_1\dots D_n\in R_0[[\epsilon]].$ We assume $D_1\dots D_n$ are linearly independent over $\bold C[[\epsilon]].$ Let  $\Cal M$ be the $\bold C[[\epsilon]]$ module generated by the $D_k.$ Then we can find an increasing sequence of integers $m_1\dots m_n$ and $E_1\dots E_n\in R_0[[\epsilon]]$ so that the $E^{[0]}_k$ are all linearly independent over $\bold C$ and so that $\epsilon^{m_k}E_k$ are a basis of $\Cal M.$  $m_1\dots m_n$ are uniquely determined and called the characteristic numbers of $D_1\dots D_n$. The $\epsilon^{m_1}E_1\dots \epsilon^{m_n}E_n$ are called a normalized basis for $\Cal M.$ Note that $E_1^{[0]},\dots,E_n^{[0]}$ are linearly independent over $\bold C.$\enddefinition
\proclaim{Corollary 3.2.7} Suppose that $T$ is one dimensional and $\{h=0\}$ is just a point $x$ and that $h$ is a parameter near $x.$ Let $g$ be the function on $\bold C^n$ defined by $$g(L)=\bold G(x,L).$$  In the context of the above definition  suppose that $(\bold G,h;\chi,\chi_1,\dots \chi_n$ is a representation of $D_1,\dots D_n$ and $g$ does not satisfy any non-zero linear combination of $\bold K_1,\dots \bold K_n, E^{[0]}_1,\dots E^{[0]}_n.$ Further assume that $$\bold v_{2j+1}(g) \equiv \sum_{i\leq j}\beta_{i,j}\bold K_i(g)\mod h$$ with $\beta_{j,j}\neq 0.$ Then the characteristic numbers of $D_1,\dots ,D_n$ are $1,3,\dots 2n-1.$ Further, the span of $E^{[0]}_1,\dots E^{[0]}_n$ is the same of the span of $\bold K_1,\dots \bold K_n.$\endproclaim
\demo{Proof} There are polynomials $P_{i,j}\in \bold C[\epsilon]$ so that 
$$\sum_{i= 1}^j P_{i,j}(h)\chi_i\equiv h^{2j-1}\bold v_{2j+1} \mod h^{2j} $$ with $P_{j,j}(0)\neq 0.$ Then let $$E'_j=\sum_{i= 1}^j P_{i,j}(\epsilon)D_i.$$  Let $\epsilon^{r_j}F_j$ be the leading term of $E_j'$ as a power series in $\epsilon.$ Then $F_j$ is a linear combination of $E^{[0]}_1,\dots E^{[0]}_n.$ If $r_j < 2j-1$, then $$F_j(g)=0.$$ So we would have to have $F_j=0.$ We conclude that $r_j\geq 2j-1.$ On the other hand, $$\bold v_{2j+1}(g) \equiv \sum_{i\leq j}\beta_{i,j}\bold K_i(g)\mod h$$ with $\beta_{j,j}\neq 0.$So $r_j\geq 2j-1$ and so $r_j=2j-1.$ Let  $$E_j=\epsilon^{-2j+1}E'_j.$$ Then $$E_j^{[0]}=\sum_{i\leq j}\beta_{i,j}\bold K_i.$$ Since the $\bold K_i$ are independent, so are the $E_j^{[0]}.$\enddemo
\subhead 3.5 \endsubhead For future reference, we will now calculate the function $\bold B$ in a very special situation.  Let $T$ be a point and let $\Cal X$ be a curve of arithmetic genus one with one node. We can find a normalization map $\pi:\bold P^1 \to \Cal X$ so that $\pi(0)=\pi(\infty)$ is the node. Let  $p$ and $s$ be  points of $\bold P^1$, $q=1/p$ and $r=1.$ We set $P=\pi(p)$, $Q=\pi(q)$, $S=\pi(s)$ and $R=\pi(r).$  Let $\delta_1$ be the image of a circle traversed counterclockwise around $0\in \bold P^1.$ We let $\Cal M =\Cal O_{\Cal X}(S).$  Then the function $\bold B$ will just depend on a number $\alpha\in \bold C$ and on $p$ and $s$. The normalized differential is just 
$$\omega_1 = \frac 1 {2\pi i} \frac {dz}{z}.$$ Explicitly, let  $$\Cal L_{\alpha}=\Cal O_{\Cal X}(\pi(x)-S).$$ Then we should have $$\frac 1 {2 \pi i} \int_S^x \frac {dz}{z}=\alpha.$$
So $$x=Se^{2\pi i \alpha}.$$ 

We use the usual parameter $z$ on $\Cal X$ coming from the parameter $z$ on $\bold P^1$ to normalize our expressions for $\lambda$ and the $s_k.$ We can now write down these expressions using a degenerate theta function
$$\vartheta_k(z) =1-\frac {z}{k}.$$ So we get 
$$\lambda(z)=\frac{\vartheta_1(z)^2 K_1}{\vartheta_p(z)\vartheta_{q}(z)},$$ where $K_1$ is chosen so that the first Laurent coefficient of $\lambda$ at $p$ is one.  
On the other hand,  $$\frac{s_{0,\alpha}}{s_{-1,\alpha}}(z)=\frac{\vartheta_p(z)
\vartheta_{x/p^2}(z)K_2}{\vartheta_q(z)\vartheta_x(z)}$$ and define $$r(z,t)= \frac
{\vartheta_{t/p^2}(z)}{\vartheta_t(z)}$$ \proclaim{Lemma 3.5.1} $$\align\bold B_{\alpha}=&\ \frac{r(q,x)r(p,s)}{r(p,x)r(q,s)}\\=&
\frac{\vartheta_{x/p^2}(q)
\vartheta_s(q)\vartheta_{x}(p)\vartheta_{s/p^2}(p)}{\vartheta_{x/p^2}(p)
\vartheta_s(p)\vartheta_{x}(q)\vartheta_{s/p^2}(q)}\\=&
\frac{(1-p/x)( 1-1/ps)(1-p /x)(1- {p^3}/{s} )    }{(1-{p^3} /x)(1-p/s)(1-1/xp)(1-p/s)}
\endalign$$\endproclaim

\subhead 3.6 \endsubhead We will now consider how to use the family $\Cal X_1 \to T_1$ to
define representations. For each positive integer $k$, we define  a
map 
$$\pi_k: \pi_*(\omega_{\Cal X_1/T_1}) \to \Cal O_{T_1} $$ by 
$$\phi_k(\omega)=Res_P(\lambda^k \omega),$$ where $Res_{P(t)}(\lambda^k
\omega)$ indicates the residue at $P(t)$ of the meromorphic section of the
dualizing sheaf of $\Cal X_{1,t}=C$ obtained by multiplying $\omega$
by $\lambda^k.$ Notice that the $\phi_k$ all vanish on $Z,$ since
$f_{-1}$ vanishes on $Z.$ We define $$\chi_k=\frac{\phi_k}{h}.$$ 
\proclaim{Lemma 3.6.1} $\psi =\chi_2+2\chi_1$ vanishes  on $Z.$
\endproclaim
\demo{Proof} This follows from Proposition 3.2.1 (5). Indeed, using the fact that $\alpha_0 \equiv -1 \mod h$, we see that $$Res_{P(t)}(\lambda_{can}^2+2\lambda_{can})
\omega $$ vanishes to order $h^2$. So $\chi_2+2\chi_1$ vanishes on $Z.$ \enddemo
  We will initially assume  $C$
is smooth.
We will also assume that  $$\int_Q^P v_k =hk.$$  we
can now associate an infinite tridiagonal matrix $C_{p}$ by defining
its $ij^{th}$ entry $C_{p,i,j}$ as
$$C_{p,i,j}=\cases 1, &\text{for $i=j+1$} \\
A_{i,p} &\text{for $i=j$} \\
B_{i,p} &\text{for $i+1=j$} \\
0 &\text{otherwise}. \endcases$$
Note that $C_{p}^n$ is a well defined infinite matrix.  If we have an
infinite matrix $E$, we define a matrix $E^+$ by 
$$E^+_{i,j}=\cases E_{i,j}, &\text{for $i<j$}\\
0 &\text{otherwise}. \endcases$$

We have identified $H^1(\Cal O_C)$ with $\bold C^n$ using the $I_k.$ So
we may think of the $\phi_k$ as sections of the map $T_1 = T \times
\bold C^n \to T.$  So we can define vertical vector fields $\chi_1,\chi_2
\ldots$ on $T_1$ corresponding to differentiating the direction
$\phi_1,\phi_2 \ldots.$  We define $\chi$ to be differentiation in the
direction $(1,2,\ldots n).$

\proclaim{Theorem 3.6.2} $$\chi_k(C_p)=[(C_p^k)^+,C_p],$$
where $[E,F]$ indicates the commutator of $E$ and $F$. \endproclaim

\demo{Proof}
Assume that $h\neq 0$ and that $h=\frac N M,$ with $N$ and $M$
integers. Then $M \gamma(Q) \in \bold Z^n$ and consequently there is a
function $\alpha$ having a pole of order $M$ at $P$ and a zero of
order $M$ at $Q.$  We normalize $\alpha$ so that $\alpha/\lambda^M$
has value one at $P$.  Then we have $$s_{k+M}=\alpha s_k$$ and
consequently, $A_k$ and $B_k$ are periodic.    Theorem 3.6.2 is just
Theorem 4 of \cite MM. On the other hand, the set of points with $h$
rational is dense and conclusion of this  theorem  holds generally.
\enddemo

We can find $T_k \in \Cal S_1 \oplus \Cal S_1$ so that 
$\Cal D_{T_k}(A_l)$ is the $l^{th}$ diagonal entry of $[(C_p^k)^+,C_p]$ and $\Cal D_{T_k}(B_l)$ is the $l^{th}$ off diagonal entry of $[(C_p^k)^+,C_p].$

\definition{Definition 3.6.3}  $D_k = D_{T_k}.$\enddefinition

For $t\in T$, let $\sigma:T\to T\times \bold C^n$ be defined by $$\sigma(t)=(t,1,2,\ldots n).$$ Then using Definition 3.2.2
$$\bold A_{k,x}=\bold f(x+kh(x)\sigma(\pi(x))$$ and 
$$\bold B_{k,x}=\bold g(x+kh(x)\sigma(\pi(x))$$
\proclaim{Theorem 3.6.4} $$(\bold f, \bold g,\chi,\chi_1\ldots \chi_n)$$ form a representation of $$D_1,\ldots D_n.$$ Further, both $\bold f$ and $\bold g$ are periodic under translation by $(1,2,\ldots n).$ Further, $$\bold f_1 =\frac{\bold f +2}{h^2}$$ and 
 $$\bold g_1 =\frac{\bold g -1}{h^2}$$ are meromorphic functions which are holomorphic at any point $t$ with $h(t)=0$ provided that the $\Cal N_t$ is non-special on the curve $\Cal X_{1,t}.$\endproclaim
\demo{Proof} We have checked the last statement when the curve $\Cal X_{1,t}$ is smooth. But if $\Cal N_t$ is non-special, then all the nearby line bundles on all the nearby curves are non-special. Hence $\bold f_1$ and $\bold g_1$ are analytic at points $w$ near $t$ provided that $h(w)\neq 0.$ Further, at points $w$ so that $h(w)=0$, the two functions are holomorphic provided that  $\Cal X_{1,t}$ is smooth. Hence both functions are holomorphic on a neighborhood of $t$, except perhaps for a subset of codimension $\ge 2.$ So by Hartog's theorem, these functions are holomorphic at $t.$\enddemo
Let $\psi = \chi_2+2\chi_1$. 
We can define $P_1$ and $P_2$ in $S_2$ so that 
$$\Cal D_{P_1,P_2}(A_k) =
A_{{k-1,p}}B_{{k-1,p}}+A_{{k,p}}B_{{k-1,p}}+2\,B_{{k-1,p}}-A_{{k,p}}B_
{{k,p}}-A_{{k+1,p}}B_{{k,p}}-2\,B_{{k,p}},
$$ 
and $$\Cal  D_{P_1,P_2}(B_k)=B_{{k,p}}\left (B_{{k-1,p}}
+A_{{k,p}}^{2}+2\,A_{{k,p}}-A_{{k+1,p}}
^{2}-2\,A_{{k+1,p}}-B_{{k+1,p}}\right ).$$
Thus we have the tame derivation $D=D_{P_1,P_2}$ of $R[[\epsilon]]$.

 Then we can use Theorem  to calculate:
$$\psi(A_{k,p})= \Cal D_{P_1,P_2}(A_k) $$ 
,
$$ \psi(B_{k,p})=\Cal  D_{P_1,P_2}(B_k).$$   

\proclaim{ Proposition 3.6.5}
 Then $\bold f,\bold g,h,\chi,\psi$ form representation of $  D_{P_1,P_2}$.
Further, $\psi$ is slow for this representation. \endproclaim

Let's calculate the first terms of $D_{\Phi}$ as a series in
$\epsilon.$
We get
$$\split
D(v^{(0)})=E_{-1}( v^{(0)})E_{-1}( w^{(0)})+
( v^{(0)})E_{-1}(w^{(0)})+\\ \quad\quad
2E_{-1}(w^{(0)})-
( v^{(0)})w^{(0)}-E_{1}( v^{(0)})w^{(0)}-2w^{(0)}.\endsplit$$
We have
$$ D_{\phi}(v^{(0)})=\Phi(D(v^{(0)}).$$
So $$\split D_{\Phi}(\epsilon^2 v^{(0)})=(E_{-1}(-2+\epsilon^2
v^{(0)})E_{-1}(1+\epsilon^2 w^{(0)})+(-2+\epsilon^2
v^{(0)})E_{-1}(1+\epsilon^2 w^{(0)})+\\ \quad\quad
2E_{-1}(1+\epsilon^2w^{(0)})-
(-2+\epsilon^2v^{(0)})(1+\epsilon^2w^{(0)})\\-E_{1}(-2+\epsilon^2
v^{(0)})(1+\epsilon^2
w^{(0)})-2(1+\epsilon^2 w^{(0)})\endsplit$$ Evaluating we get
$$ D_{\Phi}(v^{(0)}) = \epsilon(v^{(1)}- w^{(1)}) + \text{higer order
terms in }\epsilon.$$
Similarly, $$ D_{\Phi}(w^{(0)}) = \epsilon(-v^{(1)}+ w^{(1)}) + \text{higer order
terms in }\epsilon.$$

\subhead 3.7\endsubhead We will be applying this construction of slow representations to prove
Theorem 1.1. We need to have a criterion for checking that such
representations are general.  
\definition{Definition 3.7.1} Let $S$ be a complex manifold  and let $\pi:X \to S$ be a
family of stable curves parameterized by $S$. Let $s_0 \in S$ be a  point
and suppose that $N_1,N_2,\ldots N_p$ are the nodes of $X_{x_0}.$ Then
we can  find functions $f_k$ defined near $s_0$ so that the
deformation of the node $N_k$ is locally isomorphic to  $xy=f_k.$ We
say the nodes are independent at $s_0$ if  the differentials of the $f_k$ are
independent at $x_0.$ \enddefinition 
\remark{Remark}  Suppose the $g_1 \ldots g_r$ are functions on $S.$ Suppose
that the   $g_1 \ldots g_r$ have independent differentials when
restricted to some smooth $S'$ on which all the $f_k$ from the above
definition vanish. Let $S''$ be submanifold defined by the vanishing
of all the $g_i.$ Then the family $X \times_S S''$ has independent
nodes.  \endremark

We will consider the following situation: Let $\Cal M$ be a connected complex
manifold and let $\Lambda$ be a sheaf of abelian groups locally
isomorphic to $\bold Z^{2n}.$  We will assume that $\Lambda$ has a
symplectic form $\langle \ ,\ \rangle.$ Let $p$ be a point of $\Cal
M$. There is the usual
monodromy representation $\rho$ of $\pi_1(\Cal M,p)$ on the stalk
$\Lambda_p$.  Suppose that there exist $\delta_1,\delta_2,\ldots
\delta_n$ in $\Lambda_p$ so that the endomorphisms $T_i$ of
$\Lambda_p$ defined by $$T_i(\gamma)=\gamma+\langle \gamma,\delta_i
\rangle \delta_i$$ are all in the image of $\rho.$ We also assume that
$\langle \delta_i , \delta_j \rangle=0$  for all $i$ and $j$.  Let
$$\delta=\delta_1 + 2\delta_2+ \ldots +n\delta_n.$$  

 Let $\pi: V \to \Cal M$ be a vector bundle of rank $n$.  Here $V$ is
a physical bundle, i.e. $V$ is a complex manifold and for $q \in \Cal
M$, the fibers of $\pi$, $\pi^{-1}(q)=V_q$,  are given the structure of complex vector
spaces.  Let $\Cal V$ be the sheaf of analytic sections of $V$, so
that $\Cal V$ is locally free of rank $n$ on $\Cal M$.  Suppose that
$\Lambda$ is a subsheaf of $\Cal V$ so that each $\lambda \in
\Lambda_p$ can be considered a local section of $V$ defined around
$p$.  In particular, by evaluating at $p$, we get a map $\mu_p$ from
$\Lambda_p \to V_p.$ We assume that the image of $\mu_p$ is a lattice in
$V_p.$ We assume that the images of the $\delta_i$ give a complex
basis of $V_p$ for each $p.$

Let $f$ be a meromorphic function on $V.$ In particular, we can look
at $f_p$, the restriction of $f$ to $V_p.$  We will assume that $f_p$
is invariant under translation by elements of $\Lambda_p$ for all $p$.
\definition{Definition 3.7.2} Let $W$ be a complex subbundle of $V$.  We say that
$W$ is good if $f$ is constant on all the cosets of $W_p \subset V_p$
and the elements of $\Lambda_p \cap W_p$ span W as a complex vector
space. \enddefinition   

Let $U \subset \Cal M$ be an connected open set and let $\lambda$ be a section
of $\Lambda$ over $U$.  The set of $p\in U$ so that $\lambda(p) \in
W_p$ is either all of $U$ or is defined by non-trivial analytic
conditions and so is nowhere dense.  Consequently, we may find a set
of the second category $T\subset \Cal M$ so that if $p\in T $ and
$\lambda(p) \in W_p$, then $\lambda(q) \in W_q$ for all $q \in U.$ We
call such a point very general. 

If $W$ is good, then the sheaf $\Lambda_W=\Lambda \cap W \subset \Lambda$ is a
locally isomorphic to $\bold Z^k$ for some $k$ and $\Lambda_W \otimes
\Cal O_{\Cal M} =W.$ 

\proclaim{Lemma 3.7.3} Suppose that the subbundle of $V$ generated by
$\delta$ is good. Then $f_p$ is constant. \endproclaim

\demo{Proof} It suffices to prove the assertion for a very general
point $p$.  Let $W$ be a good subbundle. First note that the monodromy
representation $\rho$ leaves $L_p=\Lambda_p \cap W_p$ invariant, since $p$
is very general.

Let $Q_p$ be a maximal complex  subspace of $V_p$ so that $f_p$ is constant
on the cosets of $Q_p.$ Note that if $f_p$ is constant on the cosets of
$Q_1$ and $Q_2$, then there $f_p$ is constant on the cosets of
$Q_1+Q_2,$ so such a maximal $Q_p$ exists.  First, suppose that $L=Q_p
\cap \Lambda_p$ is not a lattice in $Q_p$.  Let $q:V_p \to
V_p/\Lambda_p$ be the quotient map and consider the closure $X$ of
$q(Q_p)$ in the torus $ V_p/\Lambda_p.$ Note that $f$ can be considered
a function on $ V_p/\Lambda_p$ which is constant on the cosets of
$q(Q_p)$ and hence on the cosets of $X.$ But $X$ is a closed subgroup of
$ V_p/\Lambda_p$ and hence there is a real subspace $Q_p'$ of $V_p$ so
that $q(Q_p')=X$ and $f_p$ is constant on the cosets of $Q_p'$. Further, $Q_p'\cap \Lambda_p$ is a lattice in $Q_p'$. $f_p$ is
meromorphic, so $f_p$ is constant on the complex subspace $Q_p''$ spanned by
the vectors in $Q_p'$.  Thus $Q_p''=Q_p$, so  $L_p=Q_p
\cap \Lambda_p$ is  a lattice in $Q_p$.  

For any point $p$, there is a
simply connected neighborhood $U\subset \Cal M$ of $p$ and a set 
$T\subset U$ of the second category in $U$
so that $Q_q\cap \Lambda_q=Q_r\cap \Lambda_r$ for all $q$ and $r$ in $T$,
since the set of subgroups of $\bold Z^{2n}$ is denumerable. Here we have identified $\Lambda_r$ with $\Lambda_q,$ since they are both identified with the global sections of $\Lambda$ over $U$. We can
then find a subsheaf $\Cal W$ of $\Cal V$ so that $\Cal W_q=Q_q$ for
all $q\in T.$ Our sheaf $\Cal W$ has been constructed in a neighborhood of
a arbitrary point $p$,  but these sheaves constructed at different
points coincide on the overlaps, since their fibers coincide over sets of the
second category.   Consequently, we have a sheaf $\Cal W$ so that the fiber
of $\Cal W$ at $q$ is $Q_q$  for a dense set of $q$. Let $\Lambda_W=\Cal W\cap \Lambda.$ Let $W\subset V$ be the physical bundle associated with $\Cal W.$  Then $f$ is constant on the cosets of $W_q$ for a dense set of $q.$
Consequently, $\Cal W_p \subseteq Q_p$ for all $p$ with equality
for a set of the second category $T$.  

Assume that $p\in T.$ Then monodromy operates on $\Cal W_p = Q_p.$  Let
$U_p$ be the real span of the $\delta_i$ .
The $\delta_i$ form a complex basis of $V_p$, so the real dimension of
$Q_q\cap U_p$ is less than  or equal to the complex dimension of $Q_p.$   Also note that
we can find a symplectic form $\langle \ ,\ \rangle$ on $V_p$ as a real vector space extending the
form on $\Lambda.$  Observe that if $v \in Q_p$ and $\langle
v,\delta_k \rangle \neq 0$, then $\delta_k \in Q_p,$ since monodromy acts. Consider the map
$T:Q_p \to Q_p \cap U_p$ defined by $$T(v) = \sum_k \langle
v,\delta_k \rangle \delta_k.$$  The observation shows that $T$ does
map to $U_p.$  Since $U_p$ is maximal isotropic, the kernel of $T$ is
contained in $Q_p \cap U_p.$ So the real dimension of $Q_p$ is less
than or equal to twice the real dimension of $Q_p \cap U_p$ and the map $T$ is
onto. The dimension of the image of $T$ is the number of $k$ so that 
$\langle
v,\delta_k \rangle \neq 0$ for some $v\in Q_p.$  Hence 
 $Q_p \cap U_p$ has a real basis consisting of some subset of the
$\delta_k$. But $\delta \in  Q_p \cap U_p$ and $\delta$ is linear
combination of the $\delta_k$ so that all the $\delta_k$ appear
non-trivially in $\delta.$ So $Q_p=V_p$ and $f_p$ is constant.

\enddemo
\head 4. Explicit Construction of Curves \endhead

\subhead 4.1 \endsubhead We will be considering families of curves in $\bold P^1 \times \bold
P^1$ over $\bold C$, which are generically double sheeted coverings of
the second factor.  Let $$T_0 = H^0(\bold P^1 \times \bold P^1, \Cal
O(n+1,2)).$$  So an element of $T_0$ is a polynomial in the variable
$X_0,X_1,Y_0,Y_1$, which is homogeneous of degree $n+1$ in $X_0,X_1$
and homogeneous of degree 2 in $Y_0,Y_1$. Usually, we will use affine
coordinates to describe the elements of $T_0$ we will be considering ,
where we set $X_0=1$, $X_1=x$, $Y_0=1$ and $Y_1=y.$ One can easily
pass from the affine coordinates to the homogenous coordinates.  So
the elements of $T_0$ of interest to us can be described by
$$a(x)y^2+b(x)y+c(x).$$ Let $T_1 \subset \bold P^1 \times \bold P^1
\times T_0$ be the universal curve and let $\pi_3:T_1 \to T_0$ and $\pi_1:T_1\to \bold P^1$ be the obvious
projections.

There is an birational involution $\iota$ of $T_1$ over $T_0$ given by 
$$\iota(x,y)=(x,-{b(x) \over a(x)}-y).$$  
Note that $\iota$ is only defined 
for those points $(x,y)$ with $a(x) \neq 0.$

Consider the map $\Lambda$ from $\bold C^n$ to $T_0$:
$$\Lambda(\alpha_1,\alpha_2,\ldots \alpha_n)=
(y^2-x)(x-\alpha_1)\ldots (x-\alpha_n)=P_{\alpha_1,\ldots \alpha_n}.$$ 
 Let $L_0$ be the image of $\Lambda$ in $T$.
Let $$C(\alpha_1,\alpha_2,\ldots \alpha_n)$$ be the curve  
$$\pi^{-1}_0(\Lambda( \alpha_1,\alpha_2,\ldots \alpha_n)).  $$ 
The $\alpha_k$ form a partial system of local coordinates around any point 
$\Lambda(\alpha_1,\alpha_2,\ldots \alpha_n)$ on $T_0$ as long as the $\alpha_k$
are distinct
and so $L_0$ is a submanifold of $T_0$ at those points.

Let $P_0$ denote $\Lambda(1, \frac
1 {2^2},\frac 1 {3^2}, \ldots \frac 1 {n^2}) \in T_0.$ We will be
investigating curves in a neighborhood of
$$C_0=C(1, \frac 1 {2^2},\frac 1 {3^2}, \ldots \frac 1 {n^2}).$$ 
For $P$ near $P_0$, $T_1\to T_0$ forms a family of semi-stable curves
parameterized by a neighborhood $U \subset T_0$ of $P_0.$ On $C_0$, we
have $2n$ nodes $N_k=(k^2,k)$ for $k$ between $-n$ and $-1$ and
between $1$ and $n.$ 

\proclaim{Lemma 4.1.1} The nodes of the family $T_1 \to T_0$ are
independent near $P_0$. .\endproclaim

\demo{Proof} By replacing $U$ by a smaller neighborhood of
$P_0$, we can find $f_k$ defined on $U$ so that the deformation of the
node $N_k$ is locally given by $$(x-\frac 1 {k^2})(y^2-x)+f_k(P)$$ for
$P\in T_0$ near $P_0.$  The $f_k$ have independent differentials at
$P_0.$ Indeed, it suffices to show that for any $k$, we can construct
a map $\psi_k:D \to T_0$, where $D$ is the unit disk, so that
$\psi_k^*(f_p)$ vanishes identically if $p \neq k$, but vanishes to
exactly order one  at $0 \in D$  if $p=k.$   $$W(x,y,t) = 
((y^2-x)(x-\frac 1{k^2})+ (t^2+\frac {2t}{k})(y+\frac 1 k)^2))\prod_{{p^2}\neq k^2} (x-\frac 1 {p^2}).$$ Note that
$$\frac{\partial W}{\partial t}(\frac 1{k^2},\frac 1 {k},0) \neq 0$$ so that the total
space family of curves over $D$ defined by $W=0$ is smooth at
$(\frac 1 {k^2},\frac 1 {k},0),$ so that  $\psi_k^*(f_k)$ vanishes exactly once at $t=0.$
On the other hand, $(\frac 1{k^2},-\frac{1}{k})$ continues to be a node of the curve
$W(x,y,t)=0$ so   $\psi_{-k}^*(f_k)=0$ and two distinct nodes continue to lie over $x=p^2$ for
$p^2\neq k^2.$ So $\psi_{p}^*(f_k)=0.$   
\enddemo

\subhead 4.2 \endsubhead For $k>0$, choose small circles $\beta_k$ oriented counterclockwise
 around the points $\frac 1 {k^2}$ and let $\delta_{k,0}$ be the lift of $\beta_k$ to
$C_0$
which is near to the point $(\frac 1{k^2},\frac 1 k)$ so that $\pi_1\circ \delta_{k,0} = \beta_k.$
For some neighborhood $U
 \subset T_0$ of $P_0,$ we can find a map 
$\delta_k: S^1 \times U \to \pi_3^{-1}(U)$ defined over $U$ which restricts to $\delta_{k,P_0},$ where $\pi_3:T_1 \to T_0$ is the projection.
Let $\omega_{T_1/T_0}$ be the sheaf of relative dualizing
 differentials. Let $\delta_{k,Q}$ be the cycle $\delta_k(S^1,Q)$ on
 $T_{1,Q}$ for $Q \in T_0.$
By possibly shrinking $U$, a section $w$ of $\omega_{T_1/T_0}$ over
 $V \subset U$ can be integrated fiberwise over the cycle $\delta_{k,Q}.$ 
 Thus we obtain a maps over $V$,
$$\int_{\delta_k}:\pi_{3,*}(\omega_{T_1/T_0})\to \Cal O_{T_0},$$ i.e. $$\left(\int_{\delta_k}\omega\right)_Q=\int_{\delta_k(S^1,Q)}\omega_Q.$$ Thus we
get a map $$\Psi:\pi_{3,*}(\omega_{T_1/T_0})\to \bigoplus_k \Cal
O_{T_0}$$ as the direct sum of the $ \int_{\delta_k}$.

Note that we can compute the dualizing differentials on
$D=C(\alpha_1,\ldots,\alpha_n).$ Namely, let $w_k$ be the differential
$$ \frac {\sqrt \alpha_k dy} {\pi i (y^2-\alpha_k)}.$$  Then $w_k$
extends to a section of $\omega_{D}$ and 
$$\int_{\delta_p}w_k=\delta_{p,k},$$ where $\delta_{p,k}$ indicates the
Kronecker delta function. So by shrinking $U$ again to a neighborhood
of $P_0$, we can assume that
$\Psi$ is an isomorphism. We have established:  \proclaim{Lemma 4.2.1}  We
can find local sections $v_p$ of 
$\pi_{3,*}(\omega_{T_1 /{T_0}})$ so that 
$$\int_{\delta_k} v_p=\delta_{k,p}.$$\endproclaim

A point $t=(E,p)$ of $T_1$ consists of an equation $E(x,y)\in T_0$ for a curve
$C \subset \bold P^1 \times \bold P^1$
and a point $p \in C.$  Let $U_1$ be a small neighborhood of $(P_0,(0,0)).$ 
We introduce functions $h_p$  of  $z \in U_1$ by the formula
$$h_p(z)=\int_{\iota(z)}^z v_p,$$
where we define this when the projection of $z$ to 
$\bold P^1 \times \bold P^1$  near
zero and $P=
\pi_3(z)$ is near $P_0$. Here 
we must specify a path $\gamma$ from $\iota(z)$ to $y$.  First we ask
that  $\pi_3(\gamma)$ be a point and that the projection of $\gamma$
to $\bold P^1 \times \bold P^1$  lie near $(0,0)$. 
 Notice that on $C_0$, the point $(0,0)$ is fixed under $\iota.$ 
So as long 
as we stay near to $(0,0)$ and $C_0$, it makes sense to ask that the path from 
stays near $(0,0)$. With these assumptions, $h_p$ is well defined. 
The functions $h_k$ all vanish on the ramification locus $\Cal R$ of
 the map $\pi_1 \times \pi_3:U \to \bold P^1 \times T_0,$ when the
 $h_k$ 
are defined, for $\Cal R$ is just defined by $\iota(z)=y.$ Also
 $\pi_3:\Cal R \to T_0$ is a local isomorphism at the points we are
 considering. 

\subhead 4.3 \endsubhead We can compute these $h_k$ on the curves $C(
 \alpha_1,\ldots,\alpha_n).$  The projection of our  path $\gamma$
 lies on the
curve  $C(
 \alpha_1,\ldots,\alpha_n).$ Let the projection of $z$ to the second
 factor of $\bold P^1 \times \bold P^1$   be $y_0.$
$$\align h_k(z)&=\int_{-y_0}^{y_0}  \frac {\sqrt \alpha_k dy} 
{\pi i (y^2-\alpha_k)}\\
&=   \frac  {\sqrt \alpha_k }  {\pi i}   
 \log(\frac {\alpha_k-y_0} {\alpha_k+y_0})\\
&\approx \frac {2y_0} {\pi i\sqrt \alpha_k}.\endalign$$
where $\approx$ indicates approximately when $y_0$ is close to zero.

This means that each of the $h_k=0$ defines $\Cal R \cap
\pi_3^{-1}(L_0)$ as a subscheme of $\pi_3^{-1}(L_0)$ in a
neighborhood of $P_0$, since $y=0$ vanishes to order one    $\Cal R \cap
\pi_3^{-1}(L_0) \subset \pi_3^{-1}(L_0).$ Since all the $h_k$ vanish
on $\Cal R$ in a neighborhood of $(0,0,y)$, we see that all the 
$h_k$ vanish to order one along $\Cal R \cap \pi_3^{-1}(L_0).$  

Let $$H_k = \frac {h_k} {h_1} $$ for $k$ from 2 to $n$.
 $\Cal R \cap \pi_3^{-1}(L_0)$ can be identified with $L_0$ locally via $\pi_3$, so we can 
use the coordinates $\alpha_1,\ldots \alpha_n$ as coordinates on 
$\Cal R \cap \pi_3^{-1}(L_0).$ Then when restricted 
to $\Cal R \cap \pi_3^{-1}(L_0)$, the $(H_k)_{\vert \Cal R\cap \pi_3^{-1}(L_0)}$ just become 
$$(H_k)_{\vert \Cal R\cap \pi_3^{-1}(L_0)}=\root\of{\frac {\alpha_k} {\alpha_1}}$$ Note that the equations
$h_1,H_2\ldots H_n$ all have independent differentials at $(0,0,P_0),$
 since the $H_k$ have independent differentials when restricted to
$\Cal R
\cap \pi_3^{-1}(L_0)$. 

\proclaim{Lemma 4.3.1} We have $$\frac{\partial}{\partial \alpha_k} H_k \neq
0$$ near $(0,0),P_0$. \endproclaim
\subhead 4.4 \endsubhead Let $$ T = \{ t \in T_1 | H_k(t) = k\} $$
and let $$\Cal X = T \times_{T_0} T_1.$$ The map $\pi:\Cal X  \to T$
has a canonical section $P: T \to \Cal X$ defined in the following
way: A point $t$ of $T$ consists of an equation $E(x,y)\in T_0$ for a curve
$C \subset \bold P^1 \times \bold P^1$
and a point $p \in C.$ $\pi^{-1}(t)$ is canonically identified with
$C.$ So we define $P(t) = p.$ We have another section $Q:T \to \Cal X$
defined by $Q(t) = \iota(p).$  On the other hand, the ramification
locus $\Cal R \subset U$ is locally isomorphic  to $T_0.$  So we can
find a section $R: U' \to \Cal X$ so that $\iota(R(t))=R(t),$ where $U'$ is a neighborhood of $(0,0,P_0).$ By
shrinking $U$, we can assume that the images of $P(t)$, $Q(t)$ and
$R(t)$ in $\bold P^1 \times \bold P^1$ are all near $(0,0).$ By localizing on $T$, we can assume that there is a section $s$ of $\Cal O(2R(t))$ which is not constant when restricted to any fiber of $\pi.$  By possibly further restricting $T$ we can a map $\bold z$ of a neighborhood of $\Cal R$ to $T\times D$ over $T$ so that 
$$s=\frac{1}{\bold z^2}.$$

\proclaim{Lemma 4.4.1} $T$ is smooth near $(P_0, (0,0))=t_0$.  $\pi:\Cal X \to
T$ is a family of semi-stable curves .  $\Cal X_{t_0}$ has $2n$ nodes
and these nodes are independent. The $I_{\delta_k}$ for $k>0$ form a basis of
the $R^1\pi_*(\Cal O)$ locally. Let $T_2$ be the subset of $T$
defined by $h=0$ and let $X_1 = \pi^{-1}(T_2).$ Let $V$ be the
physical bundle associated to $R^1\pi_*(\Cal O_{X_1})$ Let $f$ be a
meromorphic function $V$ defined on a neighborhood of the inverse
image of $t_0.$ Then if $f$ is constant on the fibers of the line
bundle generated by $$\delta =
\delta_1 + \delta_2 +\dots +\delta_n,$$ then $f$ is constant on the
fibers of $V.$
\endproclaim

At this point, we will make a choice of a line bundle $\Cal M$ on $T_1.$  Note that 
$$T_1 \subset T_0\times \bold P^1\times \bold P^1.$$ A point of $T_1$ consists of an equation $E\in T_0$ and a point $(x,y)$ with $E(x,y)=0.$ So we can map $\phi:T_1\to \bold P^1$ by $\phi(E,x,y)=y.$ On the other hand, locally around $P_0$, we can find a map $\gamma:U_0\to T_1$ so that $\gamma(U_0)\subset \Cal R.$ Thus $\gamma(U_0)$ is a divisor on the inverse image $U_1$ of $U$ in $T_1$.  We let $$\Cal M=\phi^*(\Cal O_{\bold P^1}(1))\otimes \Cal O_{U_1}(-\gamma(U_0)).$$ We denote the pullback of $\Cal M$ to $\Cal X$ by $\Cal M$ again.  We now have a function \newline$\bold B(z,\alpha_1,\alpha_2\ldots \alpha_n)$ defined locally on $T_1-T_2.$ Thus $\Cal M$ restricted to the curve $\{E=0\}=C$ is just $$\phi^*(\Cal O_{\bold P^1}(1))\otimes \Cal O_C(-R)$$ on the curve $C(\alpha_1,\alpha_2,\ldots,\alpha_n).$ This bundle has degree one on all the vertical components of $C(\alpha_1,\alpha_2,\ldots,\alpha_n),$ but degree zero on the curve $\{y^2=x\}.$

\subhead 4.5 \endsubhead Let $b$ be a non-zero  integer between $-n$ and $n$.  Let $$W_{t}(x,y) = 
((y^2-x)(x-b^2)+(t^2-2bt)(y-b)^2)).$$ Note that the point
$(b^2,b)$ is a node the curve $W_t=0$ for all $t.$  Fixing $t$, let 
$(x_1,y_1)$ be a generic point of $W_t=0.$ Let $$\Cal S=\prod_{k\neq b}(x-g_k).$$ Let  $\Cal Z=(W_t\Cal S,(x_1,y_1))\in T_1.$ Our aim is to evaluate $$\bold B(\Cal Z,\alpha_1,\alpha_2\ldots \alpha_n).$$  
We can find a map of $\psi_{1,t}:\bold P^1
\to  \bold P^1 $ by $$\psi_{1,t}(w)=-{\frac {4\,w{t}^{2}-8\,wt\,b-{b}^{2}+2\,{b}^{2}w-{b}^{2}{w}^{2}
}{\left (w+1\right )^{2}}}$$ and 
$\psi_{2,t}:\bold P^1
\to  \bold P^1 $ by $$\psi_{2,t}(w)=-{\frac {4\,w{t}^{2}-8\,wt\,b+bt-2\,{b}^{2}-b{w}^{2}t+2\,{
b}^{2}w}{\left (wt-t+2\,b\right )\left (w+1\right )}}.$$ Let
$\psi_{t}=(\psi_{1,t},\psi_{2,t}):\bold P^1 \to \bold P^1
\times \bold P^1. $  Then $\psi_{t}$ maps $\bold P^1$ to the curve
$W_t=0$ and is in fact the normalization of this curve for $t$
generic. Further, we have $\psi_{1,t}(0)=\psi_{1,t}(\infty)=b^2$
and 
$\psi_{2,t}(0)=\psi_{2,t}(\infty)=b$, so 0 and $\infty$ map to
the node of $W_t=0$.  Further,  $\psi_{1,t}$ is ramified at 1
and in fact $  \psi_{1,t}(w)=\psi_{1,t}(1/w).$
Now $\psi_{2,t}^{-1}(\infty) = \{-1,1-2b/t\}=\{P_1,P_2\}.$  
Let $C\subset \Cal X_p$ be the curve with
equation $W_t =0.$

Let $\Cal L$ be any line bundle on $\Cal X_z$ 
which has degree zero on all the irreducible components of  $\Cal X_z.$
\proclaim{Lemma 4.5.1} The natural restriction map  $\phi:H^0( \Cal X_z,\Cal
L\otimes \Cal M_z)\to H^0(C,\Cal L \otimes   \Cal M_z\otimes\Cal O_C)$
is an isomorphism.  \endproclaim
\demo{Proof} Note that $
\Cal L \otimes   \Cal M_z\otimes\Cal O_C$ has degree one on a
curve of arithmetic genus one, while $\Cal
L\otimes \Cal M_z)$ has degree $n$ on a curve of arithmetic genus $n.$
Hence, we need only show that $\phi$ is injective. A section $s$ in
the kernel of $\phi$ is a section of $\Cal
L\otimes \Cal M_z$ which vanishes on the curve $C.$ The other
components of  $\Cal X_z$ are all fibers of the projection of $\bold
P^1 \times \bold P^1$ onto the first factor.  As such, the degree of $
\Cal M_z$ on these components is one.  But these components meet $C$
in two points. So the restriction of $s$ to these components is a
section of a bundle of degree one which vanishes at two points. Hence
the section vanishes on all the vertical components, and so $s=0.$ 

\enddemo

For $v_0\in \bold P^1$, let $\psi_t(v_0)=(x_1,y_1)$. We will develop conditions on the $g_k$ so that $\Cal Z\in T.$ In fact, we will write $g_k$ as a function of $v_0$ and $t$.    We can write 
$$w_k = 
\frac 1 {2\pi i} \left(\frac{1}{z-g_k}-\frac{1}{z-1/g_k}\right)dz$$ for $k\neq 0$, while $$w_0=\frac 1 {2\pi i}\frac {dz} {z}.$$
We will consider $g_k$ close to $$\frac
{b-k}{b+k}.$$  The $w_k$ are the pullbacks of the canonical differentials on curve in $\bold P^1\times \bold P^1$ corresponding to $\Cal Z.$ Then we have 
$$\align h_k &=\int_{1/v_0}^{v_0} w_k\\&=
 \frac 1 {2\pi i} \left(\log(v_0-g_k)-\log(1/v_0-g_k)
-\log(v_0-1/g_k)+\log(1/v_0-1/g_k)\right)\\
&= \frac 1 {\pi i} \log \left(\frac{v_0-g_k} {1-g_k v_0
 }\right)\endalign$$
 Note that we have chosen the usual branch of the $\log$ so that 
$h_k$ vanishes when $v_0=1.$ We have 
$$\align h_b &=\int_{1/v_0}^{v_0} w_b\\&= \frac 1 {\pi i} \log(
v_0)\endalign$$
Next we choose $g_1\ldots g_{b-1},g_{b+1} \ldots$ so that 
$$H_k = k.$$  We do this by first choosing $g_1$ as a function of
$v_0$ and $t$
to make $h_b = b
h_1.$  Indeed, we can just take 
$$g_1=v_0\frac{1-v_0^{b-1}}{1-v_0^{b+1}}.$$ Note that $g_1$ is
analytic even when $v_0=1$ and in fact $$g_1=\frac{b-1}{b+1},$$ when
$v_0=1.$ We can find similar formulas for $g_k$ in terms of $v_0$ for
the rest of the $k$ which are not $b.$ Let $$R(t,v_0)(x,y) =W_{t}(x,y) \prod_{j \neq b} (x-\psi_{1,t}(g_k(v_0,t))).$$Recall that $$h=\frac 1 b h_b=\frac{1}{b i\pi}
\log(v_0)$$ so that $v_0=\exp(bi\pi h).$ Then we can construct a map $$\Phi_{b}:\bold C\times \bold C \to T$$ by
$$\Phi_b(h,t)=(R(t,v_0),(x_1,y_1)).$$ Note that $\Phi_b(0,0)=(P_0,(0,0)).$
 Then $\Phi_b(h,t)\in T.$ We will use Lemma 3.5.1 to compute $$\bold B(\Phi_b(h,t),\alpha_1,\ldots \alpha_n).$$  Let $s=2b/t-1.$ Note that the pullback of $\Cal M$ to curve $C=\{W_t=0\}$ is 
$\Cal O_C(\psi_t(S)),$ where $S=\pi(s).$ Let $x=se^{2 \pi i\alpha_b}.$ Then 
 $$\bold B(\Phi_b(h,t),\alpha_1,\ldots \alpha_n)={\frac {\left (-x+v_{{0}}\right )^{2}\left (v_{{0}}s-1\right )\left (-
s+{v_{{0}}}^{3}\right )}{\left (v_{{0}}x-1\right )\left (-s+v_{{0}}
\right )^{2}\left (-x+{v_{{0}}}^{3}\right )}}
.$$ Let $$\bold H=\left(\frac{-1+\bold B(\Phi_b(h,t),\alpha_1,\ldots \alpha_n)}{h^2}\right)_{h=0}.$$  Then we can compute 
$$ \bold H = 4\,{\frac {s\left ({s}^{2}{\beta}^{2}-\beta+1-{s}^{2}\beta\right ){
\pi }^{2}}{\left (s\beta-1\right )^{2}\left (s-1\right )^{2}}},
$$
where $$\beta =\frac{x}{s}=\exp(2 i\pi \alpha_b).$$ So we get 
\proclaim{Lemma 4.5.2} $$\left(\frac{d\bold H}{dt}\right)_{t=0}=  \frac {2{\pi }^{2} (\beta-1 )b}{\beta}$$\endproclaim

\subhead 4.6 \endsubhead Recall that the family  $\pi:\Cal X \to T$ depends on the integer
$n$. (The curves have bidegree $(2,n+1)$.)  Let's rename  $T$ as $V_n$ and $\Phi_b$ as $\Phi_{b,n}$ and $\Phi_b(0,0)=Q_n.$ We also define $$\phi_{b,n}(t)=\Phi_{b,n}(0,t).$$ We also denote 
$N_n = V_n \times \bold C^n.$ We claim that we can fit the $N_n$ into
a compatible family. Our first task is to construct a map $k_n: V_n
\to V_{n+1}.$ Let $$T_{0,n}= H^0(\bold P^1 \times \bold P^1, \Cal
O(n+1,2)).$$ We map $u_n: T_{0,n}  \times \bold C \to  T_{0,n+1}$ by
$$u_n(\Cal P,\beta)(x,y) = \Cal P(x,y)(x-\beta).$$ Here we will only deal with $\beta$ near $1/(n+1)^2$ so the
curve $x - \beta =0$ will meet the curve $P(x,y)=0$ transversally. Now we
see that the fiber $C_n(\Cal P)$ of $T_{1,n}\to T_{0,n}$ over $\Cal P$ is
naturally a subcurve of the fiber of $C_{n+1}(u_n(\Cal P)),$ where $T_{1,n}$ is the universal curve over $T_{0,n}$. Thus we have natural maps
$$k_n: T_{1,n}\to T_{1,n+1}.$$ Further, the
normalized differentials $v_{k,n+1}$ on the curve $C_{n+1}(\Cal P)$
restrict to the normalized differentials $v_{k,n}$ for $k = 1 \ldots
n.$ So the functions $$h_{k,n}(\Cal P,\beta)=\int_Q^P v_{k,n+1}$$ on the curve
$C_{k_n(\Cal P,\beta)}$ are independent of $\beta$ and in fact 
$$h_{k,n}(\Cal P,\beta)=\int_Q^P v_{k,n},$$ where the latter integral is
taken on the curve $\Cal P(x,y)=0.$  On the other hand, $$\frac
{\partial}{\partial \beta} \frac {h_{n+1,n+1}}{h_{1,n+1}} \neq 0$$ when $\beta$ is near
$$\frac{1} {(n+1)^2}$$ by Lemma 4.5.2. 

Now $$V_n\subset T_{1,n}$$ is the set of $(\Cal P,P)$ with $P \in \Cal C_n(\Cal P)$ so that $$\int_{\iota(P)}^P v_{k,n}=k\int_{\iota(P)}^P v_{1,n}=kh_{1,n}.$$ The functions $k_n$ map $V_n \to V_{n+1}.$ Hence we have functions $\bold f_n$ and $\bold g_n$ on $V_n\times \bold C^n$ which form a representation of $D_1 \ldots D_n\ldots.$ Further, $$\bold f_n(v,z_1,z_1\ldots z_n)=\bold f_{n+1}(k_n(v),z_1\ldots z_{n+1}).$$We let $V_n'$, $\bold f'_n$ be the extended representations. (Definition 2.6.5.)

Let $g_n'$ be the meromorphic functions defined on $(V_n'\times \bold C)\cap \{h_{1,n}\}=0$ $$g'_n(v,z)=\bold g'_{1,n}(v,z,2z,\ldots nz).$$ 
\proclaim{Lemma 4.6.1} The family $\{g'_{n}\}$ is generic. \endproclaim
\demo{Proof} Let $b$ be an non-zero integer between $-n$ and $n.$ Then Lemma 4.5.2 shows that there are maps $\phi_{b,n}$ from the disk $D\times \bold C$ to $V_n$ so that 
$$\left(\frac{\partial \phi_{b,n}^*(g_{1,n})}{\partial t}\right)(0,z)=\frac{4\pi^2}{b}(1-\exp(-2\pi i b z)).$$ We can then consider $V_n$ as a subset of $V_n'=V_n\times \bold C$ by sending $v$ to $(v,0)$ and hence we can consider $\phi_{b,n}$ as a map to $V'_n\times \bold C$.When $b=0$, we let $$\phi_{0,n}(s,z)=((Q_r,s),z).$$ Then $$\left(\frac{\partial \phi_{0,n}^*(g_{1,n}')}{\partial s}\right)(0,z)=2.$$

\enddemo
\proclaim{Theorem 4.6.2} There is a $Q\in R_0[[\epsilon]] $so that $$D_i(v^{(0)}-Q )\in \Cal I_Q.$$\endproclaim
\demo{Proof} We have constructed a representation satisfying the hypotheses of Theorem 2.9.5.\enddemo

\proclaim{Theorem 4.6.3} The characteristic numbers of $D_1,\dots ,D_n$ are $1,3,\dots 2n-1.$ Further, the span of $E^{[0]}_1,\dots E^{[0]}_n$ is the same of the span of $\bold K_1,\dots \bold K_n,$ where the $E_k^{[0]}$ are the leading terms of a normalized basis of the $\bold C[[\epsilon]]$ module  $\Cal M$ generated by the $D_k.$ (Definition 3.2.6) 
\endproclaim

\demo{Proof} Using Lemma 2.8.3 and Lemma 4.6.1, we can construct a representation of $D_1,\ldots,D_n$ satisfying the hypotheses of Corollary 3.2.7.\enddemo

 \head 5. Poisson Structures \endhead

\subhead 5.1 \endsubhead Let $$\hat R_0 =\bold C[ \ldots \hat a_{-1},\hat a_{0}, \hat a_1 \ldots
\hat b_{-1},\hat b_{0}, \hat b_1 \ldots]$$ and let
$$\hat R = R_0 [\zeta,\zeta^{-1}].$$ We say
a monomial in the $\hat a_k$ and $\hat b_l$ has weight $r$ if the sum of
the subscripts of the  $\hat a_k$ and $\hat b_l$ sum to $r$. So the
monomial $\hat a_1 \hat a_2 \hat b_{-3} $ has weight 0, as does $\zeta.$  Let $I_k
\subset R$ be
the $\bold C$ span of all the elements of weight $k$.  Let $M_N$ be
the ideal of $R_0$ generated by $$\hat a_N, \hat a_{N+1} \ldots \hat
a_{-N}, \hat a_{-N-1} \ldots
\hat b_N, \hat b_{N+1} \ldots \hat b_{-N}, \hat b_{-N-1} ,$$ i.e. a monomial is in $M_N$ if it involves $\hat a_k$ or $\hat b_k$ with $|k| \geq N$. We also
let  $M_N$ denote the induced ideal in $\hat R$.
Let $\hat I_k$ be the completion of $I_k$ with respect to subspaces
$I_k \cap M_N$ as $N \to \infty.$ Then $$\Cal F =\bigoplus_k \hat
I_k$$ is called the Fourier ring.  $\Cal F$ is naturally a graded
ring. 

\subhead 5.2 \endsubhead Suppose we are given elements $f,g$ of $\Cal S =\bold
C[z,z^{-1}]$. We define $\hat
a_n(f,g)\in
\bold C$ to
be the coefficient of $z^n$ in $f$ and $\hat b_n(f,g)$ to be the
coefficient of $z^n$ in $g$.  If $P \in \Cal F$, we can extend these
definitions to define $P(f,g) \in \bold C[\zeta,\zeta^{-1}]$.  So if $f=\sum \alpha_n z^n $ and $g=\sum \beta_n z^n$ and $P \in I_k$, then $P(f,g)$ is the result of substituting $\alpha_n$ for $\hat a_n$ and $\beta_n$ for $\hat b_n$ in $P$.  Note that $P(f,g)$ is well defined. To check
that two elements of $\Cal F$ are equal, all we have to do is to check
that they induce identical functions on $\Cal S^2.$  Also, if $P\in \Cal F[Z]$, then we can define $P(f,g)\in \bold C[Z,\zeta,\zeta^{-1}].$ We denote by
$\Cal F_0$ the analogous construction for $R_0.$  We can naturally map
$\Psi:\Cal F \to \Cal F_0[[\epsilon]] $ by sending $\zeta$ to $\exp(2\pi i \epsilon)$
considered as a formal power series in $\epsilon.$ 

Suppose we have $f\in \Cal S$ and let $N$ be a positive integer,
which we will think of as being large and let $$\zeta_N
=\exp(\frac{2\pi i}{N}).$$ Then we can define 
$$T_N( f)(n) = f(\zeta_N^n)$$ so $T_N( f):\bold Z/N \to \bold C.$ 
Note that if $N$ is sufficiently large depending on $f$ and $k$, then 
we can recover $f$ from $T_N( f),$ namely 
$$\frac 1 N \sum T_N( f)(n) \zeta_N^{-nk}$$ is the coefficient of $z^k.$

Now suppose we have a  polynomial $P \in S_1$ (see \S 2.2 for definition of $S_1$) and let
$\Cal C_N$ be the set of $\bold C$ valued functions on $\bold Z$
periodic of order $N$.  We can then define
$$\multline P(F,G)(n)=\\P(\ldots F(n-1),F(n),F(n+1)
\ldots;\ldots G(n-1),G(n),G(n+1) \ldots)\endmultline$$
Given $f \in \Cal C_N$, we define 
$$\hat f(n) = \frac 1 { N} \sum_{k \in \bold Z/N} f(k)\zeta_N^{-nk}.$$
Now suppose we are given two elements $P_1,P_2$ of $S_1.$ We can find a
continuous derivations $\Cal E_{P_1,P_2}$ of $\Cal F$ with the
property that for all $f,g\in \Cal S$ if $P_i(T_N(f),T_N(g))=h_{i,N}$, $$
\left( \Cal E_{P_1,P_2}(\hat a_n)(f,g)\right)_{\zeta =\zeta_N}=\hat
h_{1,N}(n)$$ and $$
\left( \Cal E_{P_1,P_2}(\hat b_n)(f,g)\right)_{\zeta =\zeta_N}=\hat
h_{2,N}(n)$$ for all $N$ sufficiently large
depending on $f$ and $g$ and $n.$

We can construct a series of maps $$f_n:R[[\epsilon]] \to \Cal F_1 $$
with the properties $$f_n(v^{(k)})=(2\pi in)^k\hat a_n$$ and $$f_n(w^{(k)})=(2\pi
in)^k\hat b_n$$ and$$f_n(1)=\delta_{n,0}$$ and $$f_n(FG)= \sum_{l\in \bold Z}
f_{l}(F)f_{n-l}(G).$$ Then we have
$$f_n(\partial P)=2\pi inf_n(P).$$

Suppose we are given a tame derivation $D$ of $R[[\epsilon]]$ and a
derivation $\hat D$ of $\hat R_0[[\epsilon]].$ 
\definition{Definition 5.2.1} We say $D$ and $\hat  D$ are compatible if 
 $$f_n(D(P)) =\hat  D (f_n(P)).$$\enddefinition
\proclaim{Proposition 5.2.2} Given  $D = \Cal D_{P,Q}$ for $P,Q \in S_1$,
then there is a unique 
compatible $\hat  D$.  $\hat  D$ maps the image of $\Psi$ to itself and restricts to $\Cal E_{P,Q}$ on the image of $\Psi.$ \endproclaim

\subhead 5.3 \endsubhead Suppose we have a finite collection $\Cal P$ of elements of $S_1$ $$P_{-k}\ldots
P_0 \ldots P_k ,Q_{-k}\ldots Q_0 \ldots Q_k, R_{-k}\dots R_0 \ldots
R_k.$$ Any polynomial with higher index is considered to be 0. Under some conditions on $\Cal P$, we can attempt define a
Poisson bracket $\{\ ,\ \}_{\Cal P}$ on the functions $\Cal G_N$ on
$\Cal C_N^2$ by asking that the bracket be a derivation in each slot,
be anti-symmetric and  satisfy Jacobi's identity. Further, let $A_k, B_k
\in \Cal G_N$ be defined by $$A_k(f,g)=f(k)$$ and $$B_k(f,g)=g(k).$$ 
Then we can define 
$$\{A_k,A_l\}_{\Cal P}= P_{k-l}(\ldots \hat A_k, A_{k+1} \ldots;\ldots \hat A_l
\ldots)$$
$$\{A_k,B_l\}_{\Cal P}= Q_{k-l}(\ldots \hat A_k, A_{k+1} \ldots;\ldots \hat B_l
\ldots)$$
and 
$$\{B_k,B_l\}_{\Cal P}= R_{k-l}(\ldots \hat B_k, B_{k+1} \ldots;\ldots \hat B_l
\ldots).$$ Note the $\hat {\phantom x}$ in the above equations is the place holder. We will suppose that $\{\phantom x,\phantom y\}_{\Cal P}$ defines a Poisson bracket on $\Cal G_N$ for $N$ sufficiently large. We can define a modified bracket
$\{\phantom x,\phantom y\}_{\Cal P, N}$ by 
$$\{A_k,B_l\}_{\Cal P,N}= Q_{k-l}(\ldots ,-2+\frac 1 {N^2} A_k,-2+\frac 1 {N^2}  A_{k+1} \ldots;\ldots ,1+\frac 1 {N^2} B_l,
\ldots)$$

Next define $$\hat A_{k,N} = \frac 1 {  N}\sum_{l \in \bold Z/N}
\zeta_N^{-kl} A_l$$ and $$\hat B_{k,N} = \frac 1 { N}\sum_{l \in \bold Z/N}
\zeta_N^{-kl} B_l$$ 
\proclaim{Proposition 5.3.1} Suppose that $\{{\phantom x},{\phantom x}\}_{\Cal P}$ defines a Poisson bracket on $\Cal G_N$ for $N$ sufficiently large. Then there is a Poisson bracket $\{\phantom x , \phantom y \}_{\Cal P}:\Cal F
\times \Cal F \to 
\Cal F[Z]$ so that for  $f,g \in \Cal S,$ then
$$ \left(\{\hat a_k,\hat a_l\}_{\Cal
P}(f,g)\right)_{\zeta=\zeta_N,Z=1/N}=\{\hat A_{k,N},\hat A_{l,N}\}_{\Cal
P,N}(T_N(f),T_N(g))$$ for all $N$ sufficiently large with analogous
formulas for $\{\hat a_k,\hat b_l\}_{\Cal
P}(f,g)$ and \newline$\{\hat b_k,\hat b_l\}_{\Cal
P}(f,g).$

\endproclaim
 Given $P\in S_1[Z]$, define $P_N:\Cal C_N \times \Cal C_N \to \bold C$ by $$P_N(F,G)=\frac 1 N\sum_{n \in \bold Z /N\bold Z}P(-2+\frac F {N^2},1+\frac G {N^2})(n)_{Z=1/N}.$$
\proclaim{Proposition 5.3.3} Suppose $P \in S_1[Z].$  Then there is a unique
$H_P$ in $\Cal F[Z]$ so that 
$$H_P(f,g)_{\zeta=\zeta_N,\,Z=1/N} = P_N(T_N(f),T_N(g)).$$ for $N$ sufficiently large.
\endproclaim
 \proclaim{Proposition 5.3.4} $$(\{\hat a_k,H_P\}(f,g))_{\zeta=\zeta_N,\,Z=1/N}=\{\hat A_k,  P_N\}(T_N(f),T_N(g)) $$ and $$(\{\hat b_k,H_P\}(f,g))_{\zeta=\zeta_N,\,Z=1/N}=\{\hat B_k,  P_N\}(T_N(f),T_N(g)). $$\endproclaim
 Suppose that $P\in S_1[Z]$ and $f$ and $g$ are in $\Cal S$. The function of $\epsilon$ defined by $$H(\epsilon)=P(f,g)_{\zeta=\exp(2\pi i \epsilon),\, Z=\epsilon}$$ is an analytic function of $\epsilon$ and for $N$ sufficiently large, $$H(\frac 1 N)=P_N(T_N(f),T_N(g)).$$  If $$|H(\frac 1 N)|<K(f,g)N^{-l},\tag 5.3.4.1$$ for all $N$ sufficiently large, then the first $l-1$ derivatives of $H$ vanish. We can attach a formal power series $\Psi(P)$ to $P$ in $\Cal F_0[[\epsilon]]$ by setting $Z=\epsilon$ and setting $\zeta = 1 + 2\pi i\epsilon+\dots=\exp(2\pi i \epsilon).$ If 5.3.4.1 holds for all $f,g\in\Cal S$, then the first $l-1$ derivatives of $\Psi(L)$ vanish.
\subhead 5.4 \endsubhead We next calculate two Poisson brackets.  Our first bracket is given
by:
 $$\{B_n,A_n\}_{\Cal P_1} =-2B_n$$ 
$$\{B_{n-1},A_n\}_{\Cal P_1} =2B_{n-1},$$  with all other brackets between the
$A_i$ and $B_j$ being zero, except for the obvious antisymmetric
versions of these two formulas.  We can now compute:
$$\align \{\hat A_{n,N} ,\hat B_{m,N}\}_{\Cal P_1}  =& \frac 1 {N^2} \{\sum_l
\zeta^{-ln}_NA_l,\sum_k\zeta_N^{-mk}B(k)\}\\ =& \frac 1 {N^2} 
\sum_{k,l} \zeta^{-(ln+mk)}_N\{A_l,B_k\}\\
=& \frac 1 {N^2}\sum_l \zeta_N^{-l(n+m)}2B_l - \sum_k \zeta_N^{-k(n+m)+n}2B_k \\
=& \frac 2 N \hat B_{n+m} (1-\zeta_N^n),
\endalign$$
and all the other brackets zero except as obviously required by
antisymmetry. 

For our second Poisson bracket $\{\ , \ \}_{\Cal P_2}$ we take
$$\{A_k,A_{k+1}\}_{\Cal P_2}= B_k$$ 
$$\{B_k,A_{k+1}\}_{\Cal P_2}= B_k A_{k+1}$$ 
$$\{B_k,A_{k}\}_{\Cal P_2}= -B_k A_{k}$$ 
$$\{B_k,B_{k+1}\}_{\Cal P_2}= B_k B_{k+1}.$$ 


Now we can calculate the appropriate Fourier brackets:
$$\align \{\hat B_{n,N},\hat B_{m,N}\} =& \frac 1 {N^2} \{\sum_l B_l
\zeta^{-lk} , \sum_k B_k\zeta^{-km}\}\\
=&  \frac 1 {N^2}\sum_{k,l} \{B_l, B_k\} \zeta_N^{-ln-km}\\
=& \frac 1 {N^2} \sum_k \{B_{k+1},B(k)\} \zeta_N^{-(k+1)n-km}
+\sum_k\{B_{k-1},B_k\}\zeta_N^{-(k-1)n-km}\\
=& \frac 1 {N^2} 
\sum_k B_kB_{k-1}\zeta^{-(k-1)n-km}-B_kB_{k+1}\zeta^{-(k+1)n-km}\\
 =& \frac 1 {N^2} \sum_k
B_kB_{k+1}(\zeta^{-kn-(k+1)m}-\zeta^{-(k+1)n-km})\\
=&  \frac 1 {N^2} \sum_{k,r,s} \hat B_r \hat B_s 
(\zeta^{(-kn-(k+1)m)+rk+s(k+1)}-\zeta^{-(k+1)n-km+rk+s(k+1)})\\
=& \frac 1 N \sum_{n+m=r+s}\hat B_{r}\hat B_s (\zeta^{-m+s} -\zeta^{-n+s}) 
\endalign$$
Next,
$$\align \{\hat B_{n,N},\hat A_{m,N}\} =& \frac 1 {N^2} \{\sum_l B_l
\zeta^{-lk} , \sum_k A_k\zeta^{-km}\}\\
=&  \frac 1 {N^2}\sum_{k,l} \{B_l, A_k\} \zeta_N^{-ln-km}\\
=& \frac 1 {N^2} \left(\sum_l \{B_{l},A_{l+1}\} \zeta_N^{-(ln+(l+1)m)}
+\sum_l\{B_{l},A_l\}\zeta_N^{-ln-lm}\right)\\
=& \frac 1 {N^2} \left( \sum_l B_{l}A_{l+1}\zeta_N^{-(ln+(l+1)m)}
-\sum_l B_{l}A_l\zeta_N^{-ln-lm}\right)\\=& \frac 1 {N^2} \left(
\sum_l B_l(A_{l+1}\zeta_N^{-m} -A_l)\zeta_N^{-ln-lm}\right)\\
=&  \frac 1 {N^2} \sum_{l,r,s} \hat B_r \hat A_s 
\left(\zeta_N^{-m +rl +s(l+1)-ln-lm}-\zeta_N^{ +rl +sl-ln-lm}\right)\\
=&  \frac 1 {N} \sum_{r+s=n+m}  \hat B_r \hat A_s
 \left(\zeta_N^{-m +s}-1\right)
\endalign$$

Finally, we compute
$$\align \{\hat A_{n,N},\hat A_{m,N}\} =& \frac 1 {N^2} \{\sum_l A_l
\zeta_N^{-lk} , \sum_k A_k\zeta_N^{-km}\}\\
=&  \frac 1 {N^2}\sum_{k,l} \{A_l, A_k\} \zeta_N^{-ln-km}\\
=&  \frac 1 {N^2}\left(\sum_{l} \{A_l,A_{l+1}\}\zeta_N^{-ln-(l+1)m}
 +\{A_l,A_{l-1}\}\zeta_N^{-ln-(l-1)m} \right)\\
 =&  \frac 1 {N^2}\left(\sum_{l} B_l\zeta_N^{-ln-(l+1)m}
 -B_{l-1}\zeta_N^{-ln-(l-1)m} \right)\\ 
=&  \frac 1 {N^2}\left(\sum_{l,r} \hat
B_r\zeta_N^{-ln-(l+1)m+rl}-\hat B_r\zeta_N^{-ln-(l-1)m+r(l-1)} \right)\\
=&  \frac 1 {N} \hat B_{n+m}(\zeta_N^{-m}-\zeta_N^{-n})
\endalign$$

We will now investigate the  bracket $\{\phantom x,\phantom y\}_2$ defined by  $\{\phantom x,\phantom y\}_2=\{\phantom x,\phantom y\}_{\Cal P_2}.$  
Note that if we have a continuous bracket $$\{\phantom x,\phantom y\}:\Cal F \times \Cal F\to \Cal F[Z]$$ then we have an induced bracket on $\Cal F_0[[\epsilon]]$ obtained by replacing $Z$ by $\epsilon$ and $\zeta$ by the formal power series $\exp(2\pi i \epsilon).$
By abuse of notation, we continue to call the induced bracket by the some name as the original bracket. 

Let $W_l=Y_0^l$ and consider $$X_l=\sum_{l=1}^{R}\frac {(-1)^l H_{W_l}Z^{2l-2}}{l+1}.$$ \proclaim{Proposition 5.4.1} Let $$X=\lim_{l\to \infty}\Psi(X_l).$$ Then $X$ is a Casimir for $\{\phantom x,\phantom y\}_{\Cal P_2}.$ The leading term of $X$ in $\epsilon$ is $\hat b_0.$\endproclaim
\demo{Proof} We will look at $\{\hat a_p,X_R\}=V_R$. Now $$ V_R(f,g)_{\zeta=\zeta_N,\,Z=1/N}=\sum_{r\in \bold Z/N\bold Z} (\sum_{l=0}^{R-1} (\frac{(-1)^l g(\zeta_N^r)^l}{N^{2l}})\{\hat A_p, B_r\}_{\Cal P_2,\,N}(T_N(f),T_N(g)).$$ On the other hand, $$\prod_{k=0}^N(1+\frac {B_k}{N^2})$$ is a Casimir for $\{\phantom x,\phantom y\}_{\Cal P_2,N},$ as was pointed out to me by Ali Kisisel. In particular for $N$ sufficiently large, $$\align 0=&\sum_{r=0}^{N} \frac{\{\hat A_p,B_r\}_{\Cal P_2,N}}{1+\frac {B_r}{N^2}}(T_N(f),T_N(g))\\=&\sum_{r=0}^{N} \frac{\{\hat A_p,B_r\}_{\Cal P_2,N}}{1+\frac {g(\zeta_N^r)}{N^2}}(T_N(f),T_N(g))\\=&\sum_{r=0}^{N} (\sum_{l=0}^{\infty} (\frac{(-1)^l g(\zeta_N^r)^l}{N^{2l}})\{\hat A_p,B_r\}_{\Cal P_2,N}(T_N(f),T_N(g)).\endalign$$ Thus if we fix $f$ and $g$ and   $l$, we can find a constant $K(f,g)$ so that $$|V_R(f,g)_{\zeta=\zeta_N,\,Z=1/N}|<K(f,g)N^{-l}$$ for $N$ sufficiently large. Thus the first $l-1$ derivatives of $V_R(f,g)$ vanish. Since this is true of any $f$ and $g$,  the first $l-1$ derivatives of $V_R$ vanish. Thus $\{\hat a_p,X\}=0$. A similar argument shows that $\{\hat b_p,X\}=0$, so $X$ is a Casimir.  \enddemo

Suppose that $P_1$, $P_2$ and $P$ are in $S_1[Z].$ Suppose further that $$\{A_k,P_N\}_{\Cal P,N}=P_1(\ldots, A_{k-1},\hat A_k,A_{k+1},\ldots;\ldots, B_{k-1},\hat B_k,B_{k+1},\ldots)$$ and $$\{B_k,P_N\}_{\Cal P,N}=P_2(\ldots, A_{k-1},\hat A_k,A_{k+1},\ldots;\ldots, B_{k-1},\hat B_k,B_{k+1},\ldots),$$ with the $\hat {\phantom x}$ indicating  place holder, not Fourier. 
\proclaim{Proposition 5.4.2}$\Cal E_{P_1,P_2}(\hat a_p)=\{\hat a_p,H_P\}_{\Cal P} $ with a similar formula for $\hat b_k.$\endproclaim 
\subhead 5.5 \endsubhead 
We get a series of derivations $\bold D_k$ of $\Cal F[[\epsilon]]$ compatible with $D_k$ (see Proposition 2.4.1).  The $\bold D_k$ all preserve the ideal $I_Q$ generated by $$f_n(v^{(0)}-Q)=L_n.$$ the other hand, we have two compatible Poisson brackets on $\Cal F_0\epsilon]].$ Let $$Z(n)=\exp(2\pi i n \epsilon).$$The first is defined by 
$$\{\hat a_n,\hat b_m\}_1= (\delta_{n,-m} +\epsilon^2 \hat b_{n+m})( 1 - \exp(2\pi\epsilon i n))$$ with all other terms zero except as dictated by the Poisson bracket axioms. Thus we obtain $$\{\hat a_n,\hat b_m\}_1=(-2\pi i \epsilon n)\delta_{n,-m}+\text{higher order terms in $\epsilon.$}$$ In particular, $$\{L_n,L_{-n}\}_1=\delta_{n,-m}(-4\pi \epsilon i n )+\text{higher order terms in $\epsilon.$}$$ The second is defined by 
$$\align \{\hat b_n,\hat b_{m}\}_2=& 
\frac 1 {\epsilon }\sum_{n+m=r+s} \Phi'(\hat b_r)\Phi'(\hat b_s)(\exp(2\pi\epsilon i (-m+s) -\exp(2\pi \epsilon i (-n+s))\\=& \frac 1 {\epsilon^2}\{\hat b_n,\hat b_m\}_{\Cal P_2}\endalign$$  with analogous expression for $\{\hat b_n,\hat a_m\}_2$ and $\{\hat a_n,\hat a_m\}_2$ from the Fourier expressions for the second bracket $\{\ , \ \}_{\Cal P_2}.$ First suppose that $n+m=0.$ We  can then compute 
$$\align\epsilon \{\hat b_n,\hat b_{-n}\}_2=& (\exp(2\pi\epsilon i (n) -\exp(2\pi \epsilon i (-n))(1+\epsilon^2\hat b_{0})^2+\\ &\quad \sum_{r+s=0,r\neq 0}\epsilon^4\hat b_r \hat b_s(\exp(2\pi\epsilon i (n+s) -\exp(2\pi \epsilon i (-n+s))
\endalign$$ If $n+m\neq 0$, then 
$$\align\epsilon \{\hat b_n,\hat b_{m}\}_2=&\epsilon^2\hat b_{n+m}(Z(n)-Z(-n)-Z(m)+Z(-m))+\text{higher order terms in $\epsilon$.}\endalign$$
So we get 
$$ \{\hat b_n,\hat b_{m}\}_2=\delta_{n+m,0}(4\pi i n)+\text{higher order terms in $\epsilon$} .$$ We have similar results of $\{\hat a_n,\hat b_{m}\}_2$ and $\{\hat a_n,\hat a_{m}\}_2.$ 

 Now $$L_n=\hat a_n-\hat b_n +\text{higher order terms in $\epsilon$}.$$ So $$\{L_n,L_{-n}\}_2=16\pi i \epsilon n \delta_{n,-m}+\text{higher order terms in $\epsilon.$}$$

\subhead 5.6 \endsubhead Ideally, our object would be to  define  induced brackets on $\Cal F[[\epsilon]]
/I_Q=\Cal F_0[[\epsilon]]$. We will define  brackets on a somewhat different  ring $\Cal S$. First, let $I_Q'$ be the ideal generated by the $L_q$ for $q \ne	 0$ and the Casimir $$X_1=\frac X {\epsilon^2}$$
 We can define a well defined bracket $\{\phantom x ,\phantom y \}_2 $  on $\Cal F_0[[\epsilon]]/I'_Q$ as follows: If $P\in \Cal F_0[[\epsilon]]$ and   
$\bar P\in \Cal F_0[[\epsilon]]/I'_Q$ is the image of $P$, then we define a   good  extension of $\bar P$ modulo $\epsilon^n$ to be an element  $P'\in \Cal F_0[[\epsilon]]$ so that $\{L_q,P'\}_k \in (I'_Q+\epsilon^n)$ for $q\neq 0.$ An extension is good if it is good modulo $\epsilon^n$ for all positive $n.$ It is easy to see that good extensions exist.  Suppose we have constructed a good extension $P_n$ modulo $\epsilon^n.$ Then we can try $$P_{n+1}=P_n +\epsilon^{n}\sum_{q \ne 0} W_q L_q$$ for $W_q \in \Cal F_0.$ Bracketing through by $L_r$ for $r \ne 0$ allows us to choose the $W_q$ uniquely so that $P_{n+1}$ is good modulo $\epsilon^{n+1}.$ Thus we can find a good extension of $\bar P.$ Given $\bar P, \bar P' \in \Cal F_0[[\epsilon]]/I'_Q$, we can define their bracket by taking  good extensions $P$ and $P'$ of $\bar P$ and $\bar P'$ and taking their bracket in $\Cal F[[\epsilon]]$ and then reducing modulo $I_Q'.$  This construction gives a well defined bracket $\{P,Q\}_{2}$ on $\Cal F[[\epsilon]]/I'_Q.$ 
For a given $n\neq 0$, define $$\hat \beta_{n} = \hat b_n + \left(\frac 1 2+\frac {n\epsilon i\pi}{2}-\frac {\epsilon^2\pi^2 n^2}{4}\right)L_n+\epsilon^2\sum_{k\neq 0}\left(-\frac n 8+\frac {3k}{8}\right)\hat b_{n-k}\frac {L_k}{k}$$ is the good extension of $\hat b^n$ modulo $\epsilon^3$ for the second bracket. For $n=0$, let $$\hat \beta_0=\frac 1 2 \left( \hat a_0+\hat b_0\right)+\frac 1 2 \sum_{k\neq 0} L_k \hat b_{-k}.$$ 
\proclaim{Theorem 5.5.1} $$\{\hat \beta_n,\hat \beta_m\}_2\equiv i\left(\pi (n-m)\hat \beta_{n+m}-\delta_{n,-m}\pi^3 n^3\mod \epsilon^3\right).$$ The ring $\Cal F[[\epsilon]]/I'_Q$ is generated topologically by the  images of $\hat \beta_k$ and $\epsilon.$ \endproclaim

\head 6. Convergence?\endhead
\subhead 6.1 \endsubhead Let $Q$ be the element of $R_0[[\epsilon]]$ we have constructed in Theorem 4.6.2.  As in \cite{G}, we can compute the coefficients $Q_n\in R_0$ of $Q.$ If $g$ is a periodic function analytic on $\bold R$, we can ask when the power series $$\sum_{k=0}^{\infty}Q_n(g)(z)\epsilon^n\tag 1.6.1$$ converges for $z\in \bold R.$   Suppose that $$|g^{(n)}(z)| < n!,\tag 6.1.2$$ where $g^{(n)}$ indicates the $n^{th}$ derivative of $g$. For $n\leq 23,$ I calculated $Q_n(g)$ using Maple and got a bound  $|Q_n(g)(z)|<K_n$ for $z\in \bold R$ by replacing each of the terms in $Q_n$ by the obvious estimate using (6.1.2). Here are the decimal values of $K_n.$
$$\matrix n & K_n \\
2 & .500\\
3 & .375\\
4 & .359\\
5& .312\\
6& .300\\
7& .289 \\
8&.283\\
9& .288\\ 
10& .285\\
11& .305\\
12& .312\\
13& .348\\
14& .387\\
15 & .452\\
16 & .634\\
17 & .756\\
18 & 1.70\\
19 & 1.95\\
20& 7.81\\
21& 8.46\\
22 & 53.2\\
23 & 55.2\\
\endmatrix$$
In order for (6.1.1) to converge for $\epsilon < 1/T$ all we would need is that $K_n<T^n.$ On the basis of the fact that we have constructed many functions $g$ coming from algebraic geometry
for which $6.1.1$ converges and the fact that the $K_n$ appear to be growing not too fast, I believe there should be some general convergence property of $Q.$   (Calculating the case $n=23$ used over a gigabyte of memory and took over 500 hours on a Sun Enterprise.)

\bye